\expandafter\chardef\csname pre amssym.def at\endcsname=\the\catcode`\@
\catcode`\@=11

\def\undefine#1{\let#1\undefined}
\def\newsymbol#1#2#3#4#5{\let\next@\relax
 \ifnum#2=\@ne\let\next@\msafam@\else
 \ifnum#2=\tw@\let\next@\msbfam@\fi\fi
 \mathchardef#1="#3\next@#4#5}
\def\mathhexbox@#1#2#3{\relax
 \ifmmode\mathpalette{}{\m@th\mathchar"#1#2#3}%
 \else\leavevmode\hbox{$\m@th\mathchar"#1#2#3$}\fi}
\def\hexnumber@#1{\ifcase#1 0\or 1\or 2\or 3\or 4\or 5\or 6\or 7\or 8\or
 9\or A\or B\or C\or D\or E\or F\fi}

\font\tenmsa=msam10
\font\sevenmsa=msam7
\font\fivemsa=msam5
\newfam\msafam
\textfont\msafam=\tenmsa
\scriptfont\msafam=\sevenmsa
\scriptscriptfont\msafam=\fivemsa
\edef\msafam@{\hexnumber@\msafam}
\mathchardef\dabar@"0\msafam@39
\def\dashrightarrow{\mathrel{\dabar@\dabar@\mathchar"0\msafam@4B}}
\def\dashleftarrow{\mathrel{\mathchar"0\msafam@4C\dabar@\dabar@}}

\def\ulcorner{\delimiter"4\msafam@70\msafam@70 }
\def\urcorner{\delimiter"5\msafam@71\msafam@71 }
\def\llcorner{\delimiter"4\msafam@78\msafam@78 }
\def\lrcorner{\delimiter"5\msafam@79\msafam@79 }
\def\yen{{\mathhexbox@\msafam@55 }}
\def\checkmark{{\mathhexbox@\msafam@58 }}
\def\circledR{{\mathhexbox@\msafam@72 }}
\def\maltese{{\mathhexbox@\msafam@7A }}

\font\tenmsb=msbm10
\font\sevenmsb=msbm7
\font\fivemsb=msbm5
\newfam\msbfam
\textfont\msbfam=\tenmsb
\scriptfont\msbfam=\sevenmsb
\scriptscriptfont\msbfam=\fivemsb
\edef\msbfam@{\hexnumber@\msbfam}

\catcode`\@=\csname pre amssym.def at\endcsname

\expandafter\ifx\csname pre amssym.tex at\endcsname\relax \else \endinput\fi
\expandafter\chardef\csname pre amssym.tex at\endcsname=\the\catcode`\@
\catcode`\@=11
\newsymbol\boxdot 1200
\newsymbol\boxplus 1201
\newsymbol\boxtimes 1202
\newsymbol\square 1003
\newsymbol\blacksquare 1004
\newsymbol\centerdot 1205
\newsymbol\lozenge 1006
\newsymbol\blacklozenge 1007
\newsymbol\circlearrowright 1308
\newsymbol\circlearrowleft 1309
\undefine\rightleftharpoons
\newsymbol\rightleftharpoons 130A
\newsymbol\leftrightharpoons 130B
\newsymbol\boxminus 120C
\newsymbol\Vdash 130D
\newsymbol\Vvdash 130E
\newsymbol\vDash 130F
\newsymbol\twoheadrightarrow 1310
\newsymbol\twoheadleftarrow 1311
\newsymbol\leftleftarrows 1312
\newsymbol\rightrightarrows 1313
\newsymbol\upuparrows 1314
\newsymbol\downdownarrows 1315
\newsymbol\upharpoonright 1316
 
\newsymbol\downharpoonright 1317
\newsymbol\upharpoonleft 1318
\newsymbol\downharpoonleft 1319
\newsymbol\rightarrowtail 131A
\newsymbol\leftarrowtail 131B
\newsymbol\leftrightarrows 131C
\newsymbol\rightleftarrows 131D
\newsymbol\Lsh 131E
\newsymbol\Rsh 131F
\newsymbol\rightsquigarrow 1320
\newsymbol\leftrightsquigarrow 1321
\newsymbol\looparrowleft 1322
\newsymbol\looparrowright 1323
\newsymbol\circeq 1324
\newsymbol\succsim 1325
\newsymbol\gtrsim 1326
\newsymbol\gtrapprox 1327
\newsymbol\multimap 1328
\newsymbol\therefore 1329
\newsymbol\because 132A
\newsymbol\doteqdot 132B
 
\newsymbol\triangleq 132C
\newsymbol\precsim 132D
\newsymbol\lesssim 132E
\newsymbol\lessapprox 132F
\newsymbol\eqslantless 1330
\newsymbol\eqslantgtr 1331
\newsymbol\curlyeqprec 1332
\newsymbol\curlyeqsucc 1333
\newsymbol\preccurlyeq 1334
\newsymbol\leqq 1335
\newsymbol\leqslant 1336
\newsymbol\lessgtr 1337
\newsymbol\backprime 1038
\newsymbol\risingdotseq 133A
\newsymbol\fallingdotseq 133B
\newsymbol\succcurlyeq 133C
\newsymbol\geqq 133D
\newsymbol\geqslant 133E
\newsymbol\gtrless 133F
\newsymbol\sqsubset 1340
\newsymbol\sqsupset 1341
\newsymbol\vartriangleright 1342
\newsymbol\vartriangleleft 1343
\newsymbol\trianglerighteq 1344
\newsymbol\trianglelefteq 1345
\newsymbol\bigstar 1046
\newsymbol\between 1347
\newsymbol\blacktriangledown 1048
\newsymbol\blacktriangleright 1349
\newsymbol\blacktriangleleft 134A
\newsymbol\vartriangle 134D
\newsymbol\blacktriangle 104E
\newsymbol\triangledown 104F
\newsymbol\eqcirc 1350
\newsymbol\lesseqgtr 1351
\newsymbol\gtreqless 1352
\newsymbol\lesseqqgtr 1353
\newsymbol\gtreqqless 1354
\newsymbol\Rrightarrow 1356
\newsymbol\Lleftarrow 1357
\newsymbol\veebar 1259
\newsymbol\barwedge 125A
\newsymbol\doublebarwedge 125B
\undefine\angle
\newsymbol\angle 105C
\newsymbol\measuredangle 105D
\newsymbol\sphericalangle 105E
\newsymbol\varpropto 135F
\newsymbol\smallsmile 1360
\newsymbol\smallfrown 1361
\newsymbol\Subset 1362
\newsymbol\Supset 1363
\newsymbol\Cup 1264
 
\newsymbol\Cap 1265
 
\newsymbol\curlywedge 1266
\newsymbol\curlyvee 1267
\newsymbol\leftthreetimes 1268
\newsymbol\rightthreetimes 1269
\newsymbol\subseteqq 136A
\newsymbol\supseteqq 136B
\newsymbol\bumpeq 136C
\newsymbol\Bumpeq 136D
\newsymbol\lll 136E
 
\newsymbol\ggg 136F
 
\newsymbol\circledS 1073
\newsymbol\pitchfork 1374
\newsymbol\dotplus 1275
\newsymbol\backsim 1376
\newsymbol\backsimeq 1377
\newsymbol\complement 107B
\newsymbol\intercal 127C
\newsymbol\circledcirc 127D
\newsymbol\circledast 127E
\newsymbol\circleddash 127F
\newsymbol\lvertneqq 2300
\newsymbol\gvertneqq 2301
\newsymbol\nleq 2302
\newsymbol\ngeq 2303
\newsymbol\nless 2304
\newsymbol\ngtr 2305
\newsymbol\nprec 2306
\newsymbol\nsucc 2307
\newsymbol\lneqq 2308
\newsymbol\gneqq 2309
\newsymbol\nleqslant 230A
\newsymbol\ngeqslant 230B
\newsymbol\lneq 230C
\newsymbol\gneq 230D
\newsymbol\npreceq 230E
\newsymbol\nsucceq 230F
\newsymbol\precnsim 2310
\newsymbol\succnsim 2311
\newsymbol\lnsim 2312
\newsymbol\gnsim 2313
\newsymbol\nleqq 2314
\newsymbol\ngeqq 2315
\newsymbol\precneqq 2316
\newsymbol\succneqq 2317
\newsymbol\precnapprox 2318
\newsymbol\succnapprox 2319
\newsymbol\lnapprox 231A
\newsymbol\gnapprox 231B
\newsymbol\nsim 231C
\newsymbol\ncong 231D
\newsymbol\diagup 231E
\newsymbol\diagdown 231F
\newsymbol\varsubsetneq 2320
\newsymbol\varsupsetneq 2321
\newsymbol\nsubseteqq 2322
\newsymbol\nsupseteqq 2323
\newsymbol\subsetneqq 2324
\newsymbol\supsetneqq 2325
\newsymbol\varsubsetneqq 2326
\newsymbol\varsupsetneqq 2327
\newsymbol\subsetneq 2328
\newsymbol\supsetneq 2329
\newsymbol\nsubseteq 232A
\newsymbol\nsupseteq 232B
\newsymbol\nparallel 232C
\newsymbol\nmid 232D
\newsymbol\nshortmid 232E
\newsymbol\nshortparallel 232F
\newsymbol\nvdash 2330
\newsymbol\nVdash 2331
\newsymbol\nvDash 2332
\newsymbol\nVDash 2333
\newsymbol\ntrianglerighteq 2334
\newsymbol\ntrianglelefteq 2335
\newsymbol\ntriangleleft 2336
\newsymbol\ntriangleright 2337
\newsymbol\nleftarrow 2338
\newsymbol\nrightarrow 2339
\newsymbol\nLeftarrow 233A
\newsymbol\nRightarrow 233B
\newsymbol\nLeftrightarrow 233C
\newsymbol\nleftrightarrow 233D
\newsymbol\divideontimes 223E
\newsymbol\varnothing 203F
\newsymbol\nexists 2040
\newsymbol\Finv 2060
\newsymbol\Game 2061
\newsymbol\mho 2066
\newsymbol\eth 2067
\newsymbol\eqsim 2368
\newsymbol\beth 2069
\newsymbol\gimel 206A
\newsymbol\daleth 206B
\newsymbol\lessdot 236C
\newsymbol\gtrdot 236D
\newsymbol\ltimes 226E
\newsymbol\rtimes 226F
\newsymbol\shortmid 2370
\newsymbol\shortparallel 2371
\newsymbol\smallsetminus 2272
\newsymbol\thicksim 2373
\newsymbol\thickapprox 2374
\newsymbol\approxeq 2375
\newsymbol\succapprox 2376
\newsymbol\precapprox 2377
\newsymbol\curvearrowleft 2378
\newsymbol\curvearrowright 2379
\newsymbol\digamma 207A
\newsymbol\varkappa 207B
\newsymbol\Bbbk 207C
\newsymbol\hslash 207D
\undefine\hbar
\newsymbol\hbar 207E
\newsymbol\backepsilon 237F
\catcode`\@=\csname pre amssym.tex at\endcsname

\magnification=1200
\hsize=468truept
\vsize=646truept
\voffset=-10pt
\parskip=4pt
\baselineskip=14truept
\count0=1

\dimen100=\hsize

\def\leftill#1#2#3#4{
\medskip
\line{$
\vcenter{
\hsize = #1truept \hrule\hbox{\vrule\hbox to  \hsize{\hss \vbox{\vskip#2truept
\hbox{{\copy100 \the\count105}: #3}\vskip2truept}\hss }
\vrule}\hrule}
\dimen110=\dimen100
\advance\dimen110 by -36truept
\advance\dimen110 by -#1truept
\hss \vcenter{\hsize = \dimen110
\medskip
\noindent { #4\par\medskip}}$}
\advance\count105 by 1
}
\def\rightill#1#2#3#4{
\medskip
\line{
\dimen110=\dimen100
\advance\dimen110 by -36truept
\advance\dimen110 by -#1truept
$\vcenter{\hsize = \dimen110
\medskip
\noindent { #4\par\medskip}}
\hss \vcenter{
\hsize = #1truept \hrule\hbox{\vrule\hbox to  \hsize{\hss \vbox{\vskip#2truept
\hbox{{\copy100 \the\count105}: #3}\vskip2truept}\hss }
\vrule}\hrule}
$}
\advance\count105 by 1
}
\def\midill#1#2#3{\medskip
\line{$\hss
\vcenter{
\hsize = #1truept \hrule\hbox{\vrule\hbox to  \hsize{\hss \vbox{\vskip#2truept
\hbox{{\copy100 \the\count105}: #3}\vskip2truept}\hss }
\vrule}\hrule}
\dimen110=\dimen100
\advance\dimen110 by -36truept
\advance\dimen110 by -#1truept
\hss $}
\advance\count105 by 1
}
\def\insectnum{\copy110\the\count120
\advance\count120 by 1
}

\font\ninerm=cmr9
\font\eightrm=cmr8

\font\tenrm=cmr10 at 10pt

\font\sc=cmcsc10

\def\msb{\fam\msbfam\tenmsb}

\def\bba{{\msb A}}
\def\bbb{{\msb B}}
\def\bbc{{\msb C}}

\def\bbh{{\msb H}}
\def\bbi{{\msb I}}

\def\bbo{{\msb O}}
\def\bbp{{\msb P}}
\def\bbq{{\msb Q}}
\def\bbr{{\msb R}}

\def\bbz{{\msb Z}}

\def\grD{\Delta}

\def\gre{\epsilon}

\def\grl{\lambda}

\def\la#1{\hbox to #1pc{\leftarrowfill}}
\def\ra#1{\hbox to #1pc{\rightarrowfill}}

\def\fract#1#2{\raise4pt\hbox{$ #1 \atop #2 $}}
\def\decdnar#1{\phantom{\hbox{$\scriptstyle{#1}$}}
\left\downarrow\vbox{\vskip15pt\hbox{$\scriptstyle{#1}$}}\right.}

\def\bowtie{\hbox to 1pt{\hss}\raise.66pt\hbox{$\scriptstyle{>}$}
\kern-4.9pt\triangleleft}
\def\hsmash{\triangleright\kern-4.4pt\raise.66pt\hbox{$\scriptstyle{<}$}}
\def\boxit#1{\vbox{\hrule\hbox{\vrule\kern3pt
\vbox{\kern3pt#1\kern3pt}\kern3pt\vrule}\hrule}}

\def\za{\vrule height6pt width4pt depth1pt}


\font\aa=eufm10

\def\Got#1{\hbox{\aa#1}}

\def\bfp{{\bf p}}

\def\bfu{{\bf u}}

\def\calo{{\cal O}}

\def\calh{{\cal H}}

\def\calm{{\cal M}} 
\def\caln{{\cal N}} 
\def\calo{{\cal O}}

\def\cals{{\cal S}}

\def\calz{{\cal Z}}

\def\Got#1{\hbox{\aa#1}}

\def\gsp1{{\Got s}{\Got p}(1)}

\font\svtnrm=cmr17

\font\bsc=cmcsc10 at 10truept 
  
\def\ks{1} 
\def\red{2} 
\def\low{3}
\def\hig{4}
\def\top{5} 
\def\hcx{6}

\centerline{\svtnrm 3-Sasakian Geometry, Nilpotent Orbits,}
\bigskip
\centerline{\svtnrm and Exceptional Quotients}
\vskip .5in
\centerline{\sc Charles P. Boyer~~ Krzysztof Galicki~~ and~~ 
Paolo Piccinni}
\footnote{}{\ninerm During the preparation of 
this work the first two authors 
were supported by NSF grant DMS-9970904. The third author was
supported by MURST and CNR.} 
\bigskip\bigskip
\centerline{\vbox{\hsize = 5.85truein
\baselineskip = 12.5truept
\eightrm
\noindent {\bsc Abstract}.
Using 3-Sasakian reduction techniques we obtain infinite families
of new 3-Sasakian manifolds $\scriptstyle{{\cal M}(p_1,p_2,p_3)}$ and
$\scriptstyle{{\cal M}(p_1,p_2,p_3,p_4)}$ in dimension 11 and 15 respectively. 
The metric cone on $\scriptstyle{{\cal M}(p_1,p_2,p_3)}$ is a generalization of
the Kronheimer hyperk\"ahler metric on the regular maximal
nilpotent orbit of $\scriptstyle{{\Got s}{\Got l}(3,\bbc)}$ whereas the cone on
$\scriptstyle{{\cal M}(p_1,p_2,p_3,p_4)}$ 
generalizes the hyperk\"ahler metric on
the 16-dimensional orbit of $\scriptstyle{{\Got s}{\Got o}(6,\bbc)}$. These are
first examples of 3-Sasakian metrics which are neither homogeneous nor toric.
In addition we consider some 
further $\scriptstyle{U(1)}$-reductions of 
$\scriptstyle{{\cal M}(p_1,p_2,p_3)}$.
These yield examples of non-toric 3-Sasakian orbifold metrics in dimensions 7.
As a result we obtain explicit 
families $\scriptstyle{{\cal O}(\Theta)}$ of compact self-dual positive
scalar curvature Einstein metrics with orbifold
singularities and with only one Killing vector field.
}}
\bigskip\bigskip\bigskip
\baselineskip = 10 truept
\centerline{\bf Introduction}
\bigskip
In 1990 Kronheimer showed that the co-adjoint orbits in 
the complex Lie algebra ${\Got g}^\bbc$ of both semi-simple
and nilpotent elements are hyperk\"ahler [Kr1, Kr2]. 
In particular these orbits carry Ricci-flat metrics.
Later, using Kronheimer's result, Swann showed that
the hyperk\"ahler structure on the nilpotent orbits is very
special [Sw]. Such orbits admit an
action of $\bbh^*$ with the orbit space
being a compact quaternionic K\"ahler orbifold of 
positive scalar curvature. Another way of expressing
this result is to say that the nilpotent orbits are metric
cones $C(\cals)$ on compact 3-Sasakian orbifolds [BGM1].
The mimimal nilpotent orbit is easily seen to be a metric
cone $C(G/K)$, where $G/K$ is a simply-connected
3-Sasakian homogeneous space of [BGM1].

Quaternionic geometry of the regular maximal nilpotent orbit
of ${\Got s}{\Got l}(3,\bbc)$ was investigated by
Kobak and Swann in great detail [KS1]. This orbit is 12-dimensional
and, in the language of 3-Sasakian geometry, it is a cone $\caln=C(\cals)$ 
on the 3-Sasakian orbifold $\cals=\bbz_3\backslash G_2/Sp(1)$. Here
$\cals$ is simply an $\bbz_3$ quotient of the 3-Sasakian homogeneous
space associated to the exceptional Lie group $G_2$, where
$\bbz_3$ is the center of $SU(3)\subset G_2$. This is
a typical example of a bi-quotient which has orbifold singularities
and the singular locus can be easily identified with a
homogeneous 3-Sasakian 7-manifold $SU(3)/U(1)$.

Furthermore, Kobak and Swann [KS1] show that this particular nilpotent
orbit can be obtained as a hyperk\"ahler quotient of another
nilpotent orbit by a circle action. In the language of 3-Sasakian geometry
it is simply a hyperk\"ahler reduction of the metric 
cone $\scriptstyle{C\big({SO(7)\over SO(3)\times Sp(1)}\big)}$ associated to the
action of the diagonal $U(1)\subset U(3)\subset SO(7)$. As 
$\scriptstyle{C\big({SO(7)\over SO(3)\times Sp(1)}\big)}$ 
can be realized as an $Sp(1)$ reduction
of the flat space $\bbh^7=C(S^{27})$ [G] we have the following hyperk\"ahler
(3-Sasakian [BGM1], quaternionic K\"ahler [GL]) quotients:
$$\matrix{C(S^{27})&&{\buildrel Sp(1) \over \Longrightarrow}&&
C({SO(7)\over SO(3)\times Sp(1)})&&{\buildrel U(1)\over \Longrightarrow}&&
C(\bbz_3\backslash G_2/Sp(1))\cr
&&&&&&&&\cr
\downarrow&&&&\downarrow&&&&\downarrow\cr
&&&&&&&&\cr
S^{27}&&{\buildrel Sp(1) \over \Longrightarrow}&& 
{SO(7)\over SO(3)\times Sp(1)}&&{\buildrel U(1)\over \Longrightarrow}&&
\bbz_3\backslash G_2/Sp(1)\cr
&&&&&&&&\cr 
\downarrow&&&&\downarrow&&&&\downarrow\cr
&&&&&&&&\cr
\bbh\bbp^{6}&&{\buildrel Sp(1) \over \Longrightarrow}&&  
{SO(7)\over SO(4)\times SO(3)}&&{\buildrel U(1)\over \Longrightarrow}&&
\bbz_3\backslash G_2/SO(4).\cr}\leqno(0.1)$$ 
In a later paper Kobak and Swann show that any classical nilpotent
orbit can be obtained as a hyperk\"ahler quotient of a flat space of
an appropriate dimension so that the above diagram is only an example [KS2].
The middle horizontal line of Diagram 0.1 is absent from the discussion in
[KS1]as the importance of 3-Sasakian geometry in this context was realized
later [BGM1]. 

The starting point of this paper is to revisit the Kobak-Swann construction
in the context of the associated 3-Sasakian geometry. We are going to
examine the quotient construction of the orbifold fibration
$\bbz_3\backslash G_2/Sp(1)\rightarrow \bbz_3\backslash G_2/SO(4)$
showing that it admits interesting generalizations. More precisely,
the Kobak-Swann quotient can be ``deformed" by introducing
weights in much the same way toric 3-Sasakian manifolds with
$b_2=1$ can all be obtained by ``deforming" the classical
homogeneous example of the fibration ${\cal S}({\bf 1})=
SU(n)/S(U(n-2)\times U(1))
\rightarrow {\rm Gr}_2(\bbc^n)$ [BGM1]. On the other hand 
our new quotients are quite different from the construction
of $\cals(\bfp)$ considered in [BGM1] as they cannot be interpreted
as bi-quotients. In Section 2 we prove

\noindent{\sc Theorem A}: \tensl
Let  ${\bf p}=(p_1,p_2,p_3)\in\bbz^3$. For any such non-zero
$\bfp$ one can define an isometric 
action of $Sp(1)\times U(1)_{{\bf p}}\subset Sp(7)$ on
$S^{27}$ with the following property: If
$0<p_1<p_2<p_3$ are pairwise relatively prime
and ${\rm gcd}(p_1\pm p_2, p_1\pm p_3)=1$ 
then the reductions ${\cal M}({\bf p})\rightarrow {\cal Z}({\bf p})
\rightarrow {\cal O}({\bf p})$ of 
$S^{27}\rightarrow \bbc\bbp^{13}\rightarrow
\bbh\bbp^{6}$ give a compact smooth
11-dimensional 3-Sasakian manifold ${\cal M}({\bf p})$ together with the
(orbifold) leaf spaces of its 
fundamental foliations ${\cal Z}({\bf p})$ and ${\cal O}({\bf p})$.
Furthermore, the reduced space 
${\cal M}(1,1,1)$ is a 3-Sasakian orbifold $\bbz_3\backslash
G_2/Sp(1)$.
\tenrm

Analysis of the symmetry structure of all the quotients together
with the associated foliations gives

\noindent{\sc Theorem B:}
The manifold ${\cal M}({\bf p})$ of Theorem A is not toric. The corresponding
leaf spaces ${\cal Z}({\bf p})$ and
${\cal O}({\bf p})$ are compact Riemannian
orbifolds with inhomogeneous Einstein metrics of positive scalar curvature.
\tenrm

In Section 3 we investigate whether ${\cal M}(p_1,p_2,p_3)$
admits further reduction by an isometric circle action. More generally
we consider $Sp(1)\times S^1\times S^1$ actions on the 27-sphere
and ask if one can get any smooth quotients. Surprisingly, no smooth examples
can be found but orbifold quotients exist in profusion. These are
interesting, since, due to 
Theorem B, they are necessarily non-toric
(more precisely of cohomogeneity 3). 
More importantly, they yield new explicit
self-dual Einstein metrics of positive scalar curvature with
only orbifold singularities and with one-dimensional group of
isometries. We get

\noindent{\sc Theorem C}: \tensl Let $\Theta\in{\cal M}_{2\times3}(\bbz)$
be any integral  $2\times3$ matrix such that each of its three
$2\times2$ minor determinants does not vanish. In addition
suppose that the sum of the all minor determinants is nonvanishing, and none of
them is equal to the some of the other two. For any such $\Theta$
there exists
a compact 4-dimensional orbifold ${\cal O}(\Theta)$ which admits
a self-dual Einstein metric of positive scalar curvature 
with a one-dimensional group of isometries. Moreover, this metric
can be constructed explicitly as a  quaternionic K\"ahler
reduction of the real Grassmannian $Gr_4(\bbr^7)$ by an isometric action
of the 2-torus $T^2_{\Theta}$ defined by $\Theta.$
\tenrm

The first examples of positive self-dual Einstein metrics on orbifolds were
obtained in [GL]. Later such metrics were considered by
Hitchin [Hi1, Hi2]. Hitchin's examples have large group of
isometries. More recently many new orbifold metrics with 
$T^2$-symmetry group were constructed in [BGMR].
The examples presented here are perhaps the first self-dual Einstein metrics
with only one Killing vector field. We are not aware of any self-dual Einstein
metrics which have only discrete isometries.

In Section 4 we consider the obvious
higher-dimensional extension of the problem. Not surprisingly, once again
the new examples involve hyperk\"ahler geometry of a
nilpotent variety. This time it is the 16-dimensional 
nilpotent orbit of ${\Got s}{\Got o}(6,\bbc)$ which is the 
Swann bundle over the quaternionic K\"ahler 
orbifold $Gr_4(\bbr^7)/\bbz_2$. In [KS2] it is shown
how all classical nilpotent orbits can be obtained as hyperk\"ahler
reductions from flat spaces (typically in more than one way). This
orbit can be constructed  as a $Sp(1)\times U(1)$
reduction of $\bbh^8$ [KS3] and thus it appears as part of a diagram
similar to the one in (0.1). In this context our construction is a systematic 
study of the general $U(1)$ actions which, at the 3-Sasakian level, produce 
smooth metrics. We prove
the following analogue of the Theorem A:

\noindent{\sc Theorem D}: \tensl
Let  ${\bf p}=(p_1,p_2,p_3,p_4)\in\bbz^4$. For any such non-zero
$\bfp$ one can define an action of $Sp(1)\times U(1)_{{\bf p}}\subset Sp(8)$ on
$S^{31}$ with the following property: If
$0\leq p_1<p_2<p_3<p_4$, any triple
$p_i<p_j<p_k$ satisfies 
${\rm gcd}(p_i,p_j,p_k)=1$, and
${\rm gcd}(p_i\pm p_j, p_i\pm p_k)=1$ then
the the reductions ${\cal M}({\bf p})\rightarrow {\cal Z}({\bf p})
\rightarrow {\cal O}({\bf p})$ of 
$S^{31}\rightarrow \bbc\bbp^{15}\rightarrow
\bbh\bbp^{7}$ give a compact smooth 
15-dimensional 3-Sasakian manifold ${\cal M}({\bf p})$ together with the
(orbifold) leaf spaces of its
fundamental foliations ${\cal Z}({\bf p})$ and ${\cal O}({\bf p})$.
Furthermore, the reduced space ${\cal M}(1,1,1,1)$ is a 3-Sasakian orbifold
$\bbz_2\backslash {\rm Spin}(7)/{\rm Spin}(4)$.
\tenrm

These 15-dimensional quotients are ``deformations" of the standard homogeneous
3-Sasakian structure on $SO(7)/SO(3)\times Sp(1)$ which projects
to the Wolf space ${\rm Gr}_4(\bbr^7)$ in the quaternionic
K\"ahler base. In higher dimensions the orbifold
${\rm Gr}_4(\bbr^{n+3})/\bbz_2$ can equally be obtained as a quaternionic
K\"ahler reduction of $\bbh\bbp^{n+4}$ by $U(1)\times Sp(1)$.
Our construction shows that the $n=4$ case, from the standpoint
of 3-Sasakian geometry, is somewhat exceptional: In higher dimension
we can only get orbifold metrics. 
The orbifold bundles ${\cal M}^{15}({\bf p})\rightarrow
{\cal Z}^{15}({\bf p})\rightarrow {\cal O}^{15}({\bf p})$
can be viewed as singular analogues of $SO(7)/SO(3)\times Sp(1)
\rightarrow SO(7)/SO(3)\times U(2)\rightarrow {\rm Gr}_4(\bbr^7).$
This makes use of the well-known isometry between ${\rm Gr}_4(\bbr^7)$ and
the space ${\rm Spin}(7)/ (Sp(1)\times Sp(1)\times Sp(1))/\bbz_2$
of the Cayley 4-planes in $\bbr^8$.

Again, standard analysis of the symmetry structure of all the quotients together
with the associated foliations gives

\noindent{\sc Theorem B$^\prime$:}
The manifold ${\cal M}({\bf p})$ of Theorem D is not toric. The corresponding
leaf spaces ${\cal Z}({\bf p})$ and 
${\cal O}({\bf p})$ are compact Riemannian
orbifolds with inhomogeneous Einstein metrics of positive scalar curvature. 
\tenrm

In Section 5 we give what topological information is available to us. In
particular, we show that as long as $\bfp$ satifies the conditions of Theorems
A and D (actually this hypothesis can be weakened) the rational cohomology of
the corresponding $\calm(\bfp)$ is independent of $\bfp.$ Finally, in section
6 we briefly mention the construction of hypercomples structures on circle
bundles over our new 3-Sasakian manifolds (orbifolds).

\noindent{\sc Acknowledgements}: The second named author would like to
thank Universit\`a di Roma ``La Sapienza", C.N.R, M.P.I-Bonn, 
and I.H.E.S-Bures sur Yvette for hospitality and support. Parts of this
paper were written during his visits there. We would also like to thank Mike
Buchner and  Andrew Swann for comments and discussion.

\bigskip
\centerline {\bf \ks. Quotient construction of 3-Sasakian structure on
$\bbz_3\backslash G_2/Sp(1)$}
\medskip
There are two homogeneous Sasakian-Einstein geometries that
are naturally associated with the exceptional Lie group
$G_2$ [BG1, BG2]. They both come from the classical Lie group
isomorphism between
$SO(4)\subset G_2$ and $Sp(1)_-\cdot Sp(1)_+\subset G_2$.
The two $Sp(1)_\pm$ subgroups are very different. One of them
has index 1 in $G_2$ and the other one has index 3. Consequently,
the quotients are not of the same homotopy type as can be seen
from the exact sequence in homotopy for the fibration
$$Sp(1)_\pm\longrightarrow G_2 \longrightarrow G_2/Sp(1)_\pm.$$
In particular, the two spaces can be
distinguished by their third homotopy groups
being trivial in one case and $\bbz_3$
in the other. One of these quotients,
which we shall denote by $G_2/Sp(1)_-$ is diffeomorphic to the
real Stiefel manifold $V_{7,2}(\bbr)=SO(7)/SO(5)$ of 2-frames in $\bbr^7$ [HL].
As $V_{7,2}(\bbr)$ is 4-connected $\pi_3(G_2/Sp(1)_-)=0$. The other
quotient denoted here by $G_2/Sp(1)_+$ is one of the 11-dimensional
3-Sasakian homogeneous spaces and $\pi_3(G_2/Sp(1)_+)=\bbz_3$.
$G_2/Sp(1)_+$ fibers as a circle bundle over a generalized
flag $\calz=G_2/U(2)_+$, which in turn is well-known to be the
twistor space of the exceptional 8-dimensional Wolf space $G_2/SO(4)$.
The second homogeneous Sasakian-Einstein manifold is a circle bundle over
the complex flag $G_2/U(2)_-$ which can be identified with the complex
quadric in the 6-dimensional complex projective space $\bbc\bbp^6$ or,
equivalently, the real Grassmannian ${\rm Gr}_2(\bbr^7)=SO(7)/SO(2)\times SO(5)$
of oriented 2-planes in $\bbr^7$. Both 
Sasakian-Einstein metrics are well-known and  have been
studied in the context of homogeneous
Einstein geometries [BG1]. We have the following diagram of Riemannian
submersions:
$$\matrix{&&&G_2&&&\cr
          &&&&&&\cr
          &&\swarrow&&\searrow&&\cr
          &&&&&&\cr
          &{G_2\over Sp(1)_+}&&&&{G_2\over Sp(1)_-}\simeq V_{7,2}(\bbr)&\cr
          &&&&&&\cr
          &\downarrow&&&&\downarrow&\cr
          &&&&&&\cr
          &\calz={G_2\over U(2)_+}&&&&{G_2\over U(2)_-}\simeq 
          {\rm Gr}_2(\bbr^7)&\cr
          &&&&&&\cr
          &&\searrow&&\swarrow&&\cr
          &&&&&&\cr
          &&&{G_2\over SO(4)}&&&\cr}\leqno{\ks.1}$$

Poon and Salamon [PS] proved that $G_2/SO(4)$ is one of the three possible
models of positive quaternionic K\"ahler manifolds in dimension 8. Later
the geometry of $G_2/SO(4)$ was examined by Kobak and Swann [KS1] who proved
the following remarkable theorem:

\noindent{\sc Theorem \ks.2}: \tensl
The quaternionic K\"ahler manifold ${\rm Gr}_4(\bbr^7)$
admits an action of $U(1)$ such that the quaternionic K\"ahler
quotient is a compact quaternionic K\"ahler orbifold
$\calo=\calo_r\cup\bbc\bbp(2)=G_2/(SO(4)\times\bbz_3)$.
\tenrm

We will first re-examine the Kobak-Swann construction from the point of
view of the 3-Sasakian geometry of a certain $SO(3)$ V-bundle over $\calo$ (or,
equivalently, the hyperk\"ahler geometry of the regular nilpotent orbit of
${\Got s}{\Got l}(3,\bbc)$). In particular, we have the following:

\noindent{\sc Theorem \ks.3}: \tensl
The 3-Sasakian homogeneous manifold $SO(7)/SO(3)\times Sp(1)$
admits an action of $U(1)$ such that the 3-Sasakian
quotient is a compact 3-Sasakian orbifold
$\calm=\calm_r\cup SU(3)/U(1)=\bbz_3\backslash G_2/ Sp(1)_+$.
\tenrm

Theorem \ks.3 is a straightforward translation of
Theorem \ks.2 into the language of 3-Sasakian geometry and we could leave
it at that. However, we 
will outline a constructive proof of this result as our description of
the corresponding quotient differs slightly from the one given in [KS1].

One can think of the homogeneous 3-Sasakian manifold
$SO(7)/SO(3)\times Sp(1)$ as the 3-Sasakian reduction
$S^{4n-1}/\!\!/\!\!/Sp(1)$ as follows [G, BGM1]:
Let $\bfu=(u_1,...,u_7)\in S^{27}$. Consider the $Sp(1)$ action given
by multiplication by unit quaternion $\lambda\in Sp(1)$ on the left
that is
$$\varphi_\lambda(\bfu)=\lambda\bfu.\leqno{\ks.4}$$
In the $\{i,j,k\} $ basis the 3-Sasakian moment
maps for this action read:
$$\mu_i(\bfu)=\sum_{\alpha=1}^7\overline u_\alpha iu_\alpha,\qquad
\mu_j(\bfu)=\sum_{\alpha=1}^7\overline u_\alpha ju_\alpha,\qquad
\mu_k(\bfu)=\sum_{\alpha=1}^7\overline u_\alpha ku_\alpha.
\leqno{\ks.5}$$
Then, the common zero-locus of the moment maps
$$N=\{\bfu\in S^{4n-1}:\ \ \ \mu_i(\bfu)=
\mu_j(\bfu)=\mu_k(\bfu)=0\}\leqno{\ks.6}$$
is the Stiefel manifold $N\simeq SO(7)/SO(3)=V_{7,4}(\bbr)$ 
of the orthonormal 4-frames in $\bbr^7$
and the corresponding 3-Sasakian
quotient $\cals=N/Sp(1)$ is Konishi's $\bbr\bbp^3$-bundle
over the real Grassmannian of
oriented 4-planes in $\bbr^7$.
We can combine Theorem \ks.3 with this description to get

\noindent{\sc Corollary \ks.7}: \tensl 
The 3-Sasakian sphere $S^{27}$
admits an action of $U(1)\times Sp(1)$ such that the 3-Sasakian
quotient is a compact 3-Sasakian orbifold 
$\calm=\calm_r\cup SU(3)/U(1)=\bbz_3\backslash G_2/ Sp(1)_+$. 
\tenrm   

We now turn to the explicit description of the $U(1)$ quotient.
Consider the following subgroups of the group of 3-Sasakian
isometries of the 27-sphere:
$$Sp(7)\supset SO(7)\supset 1\times SO(6)\supset 1\times U(3) 
,\leqno{\ks.8}$$
where $U(1)\subset U(3)$ is the central subgroup. Explicitly, we shall write
$f:[0,2\pi)\rightarrow SO(7)$
$$f(t)=\pmatrix{1&0&0&0\cr 0&A(t)&0&0\cr 0&0&A(t)&0\cr
0&0&0&A(t)\cr}\in SO(7),\leqno{\ks.9}$$
where
$$A(t)=\pmatrix{\cos t&\sin t\cr-\sin t&\cos t\cr}\leqno{\ks.10}$$
are the real rotations in $\bbr^2$. The homomorphism $f(t)$
yields a circle action on $S^{27}$ or,
equivalently, (after performing the ``$Sp(1)$-reduction" first) on 
the homogeneous 3-Sasakian manifold
$SO(7)/SO(3)\times Sp(1)$ via left multiplication
$f(t){\bfu}$ and the associated 3-Sasakian moment map can be written as
$$\nu(\bfu)=\sum_{\alpha=1,2,3}(
\overline u_{2\alpha}u_{2\alpha+1}-\overline u_{2\alpha+1}u_{2\alpha}).
\leqno{\ks.11}$$
Note here that $\nu(\bfu)$ does not depend on the $u_1$
quaternionic coordinate. 

\noindent{\sc Definition \ks.12}: \tensl
Let us define the zero level set of
this new moment map intersected with $N$, that is
$$N_\nu\equiv N\cap \nu^{-1}(0)\leqno{\ks.12}$$
\tenrm

First we observe, following Kobak and Swann [KS1] that

\noindent{\sc Lemma \ks.13}: \tensl The manifold $N_\nu$ can be identified
with $U(1)\cdot G_2 =(S^1\times G_2)/\bbz_3$ where $U(1)\cap G_2=\bbz_3.$
\tenrm

\noindent{\sc Proof}: The argument is similar here to the one
used by Kobak and Swann in [KS1] and it is based on the Proposition 1.10
of [HL]. First, using the basis $\{i,j,k\}$ of unit imaginary quaternions,
we write $u_\alpha=u^0_\alpha+
iu^1_\alpha+ju^2_\alpha+ku^3_\alpha$ and introduce the $4\times7$ real matrix
$$\bba=\pmatrix{u_1^0&u_2^0&u_3^0&u_4^0&u_5^0&u_6^0&u_7^0\cr
u_1^1&u_2^1&u_3^1&u_4^1&u_5^1&u_6^1&u_7^1\cr
u_1^2&u_2^2&u_3^2&u_4^2&u_5^2&u_6^2&u_7^2\cr
u_1^3&u_2^3&u_3^3&u_4^3&u_5^3&u_6^3&u_7^3\cr}\equiv\pmatrix{f^0\cr f^1\cr
f^2\cr f^3\cr},\leqno{\ks.14}$$
where to make the connection with the notation in [KS1]
we also think of the rows of $\bba$ as purely imaginary
octonions ${\rm Im}(\bbo).$ In the standard basis of ${\rm Im}(\bbo)$ 
we write $f^a= u_1^ai+u_2^aj+u_3^ak+u_4^ae+u_5^aie+u_6^aje+u_7^ake.$
Let $\phi(a,b,c) =<ab,c>$ denote the 3-form defining the associative
calibration [HL] on ${\rm Im}(\bbo)$ where $<\cdot,\cdot>$ denotes the
standard Euclidean inner product. Then writing $\nu=\nu_1i +\nu_2j+\nu_3k$ a
straightforward computation shows that for $a=1,2,3$ and $\gre^{abc}$ the
totally antisymmetric tensor satisfying $\gre^{123}=1$
$$\nu_a=2<f^0f^a+\gre^{abc}f^bf^c,i>, \leqno{\ks.15}$$
where the summation convention on repeated indices is used.
Now, $\bfu\in N$ if and only if  the rows
$\{f^0,f^1,f^2,f^3\}$ of $\bba$ 
form an orthonormal frame in $\bbr^7\simeq {\rm Im}(\bbo),$ and              
one can identify $G_2$ with a special kind of oriented orthonormal
4-frame, namely those which are co-associative. This means that the 3-plane
that is orthogonal to the 4-plane defined by the frame $\{f^0,f^1,f^2,f^3\}$
is spanned by an associative subalgebra of ${\rm Im}(\bbo).$  Then one shows
that these special 4-frames satisfy the $U(1)$-moment map equation
$\nu(\bfu)=0$ and, hence, $U(1)\cdot G_2\subset N\cap
\nu^{-1}(0)$. As, $U(1)\cap G_2=\bbz_3$ it is enough to show that by acting
with $U(1)$ one gets the whole $N\cap \nu^{-1}(0)$. The argument is similar
to the one presented in [KS1]. (See [KS1] Lemma 5.1 and the discussion 
that follows.) \hfill\za

Now, Theorems \ks.2 and \ks.3 and Corollary \ks.7 
all follow from the above lemma
as we get the quotient
$$\calm={N_\nu\over U(1)\times Sp(1)}\simeq{U(1)\cdot 
G_2\over U(1)\times Sp(1)}\simeq
\bbz_3\backslash G_2/Sp(1).\leqno{\ks.16}$$

\noindent{\sc remark \ks.17}: The $U(1)\times Sp(1)$ action on the level
set $N_\nu$ is not locally free. If we divide by $Sp(1)$ first
and consider the $U(1)$ action on the orbit space $N_\nu/Sp(1)$ this
circle action is quasi-free. This means that there are only two kinds 
orbits: regular orbits with the trivial isotropy group and singular orbits
(points) where the isotropy group is the whole $U(1)$. In such cases
the quotient space is often an orbifold (or even a smooth manifold).
The stratification of the Theorem 2 is precisely with respect to the
orbit types as will be seen in the next section.

\bigskip
\def\cals{{\cal S}}

\centerline {\bf \red. Generalizations of the Kobak-Swann Quotient}
\medskip
In this section we will consider the simplest possible family of
quotients which generalize the construction of Kobak and Swann
via an introduction of weights. Instead of considering the
central $U(1)\subset U(3)$ in \ks.9 we can consider
and arbitrary circle subgroup of the maximal 
torus $U(1)\subset T^3\subset U(3)$.
Again, to be more specific, we have the following inclusions:
$$Sp(7)\supset SO(7)\supset 1\times SO(6)\supset 1\times SO(2)\times
SO(2)\times SO(2).\leqno{\red.1}$$
We can consider arbitrary circle subgroups of the last 3-torus. Let
${\bf p}=(p_1,p_2,p_3)\in \bbz^3$ and define the following homomorphism
$$f_{{\bf p}}(t)=\pmatrix{1&0&0&0\cr 0&A(p_1t)&0&0\cr 0&0&A(p_2t)&0\cr
0&0&0&A(p_3t)\cr}\in SO(7),\leqno{\red.2}$$
where
$$A(p_it)=\pmatrix{\cos(p_it)&\sin(p_it)\cr-\sin(p_it)&\cos(p_it)\cr}
\in SO(2),\ \ \ \ \ i=1,2,3 \leqno{\red.3}$$
are 2-dimensional real rotations. Note that for $p_1=p_2=p_3=1$ we recover the
example of the previous section. The homomorphism $f_{{\bf p}}(t)$
yields a circle action on the homogeneous 3-Sasakian manifold
$SO(7)/SO(3)\times Sp(1)$ via left multiplication
$f_{{\bf p}}(t){\bf u}$ and the moment map can now be written as
$$\nu_{{\bf p}}(u_1,...,u_7)=\sum_{\alpha=1,2,3}p_\alpha(
\overline u_{2\alpha}u_{2\alpha+1}-\overline u_{2\alpha+1}u_{2\alpha}).
\leqno{\red.4}$$

\noindent
Observe that without loss of generality we can assume all weights
to be non-negative as $p_i$ can be changed to $-p_i$ by renaming
the quaternions in the associated pair $(u_{2i},u_{2i+1})$.
We begin analysis of this quotient by considering the level set 
of the moment map 

\noindent{\sc Definition \red.5}: \tensl Let $N_\nu(\bfp)\subset S^{27}$
be the level set of the 3-Sasakian moment map of the 
$Sp(1)\times U(1)_\bfp$-action, {\it i.e.,} 
$N_\nu({\bf p})\equiv N\cap \{\nu_{{\bf p}}^{-1}(0)\}.$
\tenrm

We want to consider a stratification of the level set $N_\nu(\bfp)$ that will
allow us to analyze the quotient space
$${\cal M}({\bf p})= 
{N_\nu({\bf p})\equiv V_{7,4}(\bbr)\cap \{\nu_{{\bf p}}^{-1}(0)\}\over  
Sp(1)\times U(1)_{{\bf p}}}. \leqno{\red.6}$$ 
Since $N_\nu({\bf p})$ is a submanifold of the Stiefel manifold $V_{7,4}(\bbr)$
at most 3 quaternionic coordinates can vanish on $N_\nu({\bf p}).$ So setting
various quaternionic coordinates equal to zero determines a stratification of 
$N_\nu({\bf p})$ in which the strata of minimal dimension play an important
role.  We call these strata {\it vertices} although, as we shall see, they each
have two connected components. 

\noindent{\sc Lemma} \red.7: \tensl Let $0<p_1<p_2<p_3$. At a
vertex neither $u_1$ nor any of the three pairs of quaternions
$(u_{2i},u_{2i+1})$, $i=1,2,3$ can vanish. Thus, there are precisely eight
vertices and they are all diffeomorphic to $O(4).$ \tenrm 

\noindent{\sc Proof}: Every vertex has precisely 3 quaternionic coordinates
vanishing, so the Stiefel manifold becomes $V_{4,4}(\bbr)=O(4).$ Let $V$ be a
vertex.
Then $V=V_{ijkl}$ can be
represented by a matrix of the form
$$\bbb=\pmatrix{u_i^0&u_i^1&u_i^2&u_i^3\cr 
u_j^0&u_j^1&u_j^2&u_j^3\cr 
u_k^0&u_k^1&u_k^2&u_k^3\cr 
u_l^0&u_l^1&u_l^2&u_l^3\cr}, \qquad \bbb^t\bbb=\bbb\bbb^t={1\over4}\bbi_4, 
\leqno{\red.8}$$ 
where $1\leq i<j<k<l\leq7$. Here we have written a quaternionic
coordinate as $u_i=u_i^0+iu_i^1+ju_i^2+ku_i^3.$
Suppose $i>1$, {\it i.e.,} $u_1=0$ on
on $V_{ijkl}$. Then there are two possibilities: (1) the quadruple
$i<j<k<l$ contains only one quaternionic pair $(2\alpha, 2\alpha+1)$
or (2) it contains two such pairs. Let us examine the first possibility.
Without loss of generality we can take $i=2,j=3,k=4,l=6$.
With this choice the 
vanishing of the $U(1)$ moment map \red.4 now becomes 
$${\rm Im}(\bar{u}_2u_3)=0.$$
One can easily check that orthogonality of the basis forces 
$${\rm Re}(\bar{u}_2u_3)=u_2^1u_3^1+ 
u_2^2u_3^2+u_2^3u_3^3+u_2^4u_3^4=0.$$ 
But this implies $\bar{u}_2u_3=0$ which forces either $u_2$ or $u_3$ to
vanish, giving a contradiction. Now assume the second
case, that is, that $i<j<k<l$ consists of two
quaternionic pairs $(2\alpha, 2\alpha+1)$. Again, when we take
$i=2,j=3,k=4,l=5$, the orthogonality of the vectors in the associated
$\bbb$-matrix forces
$$\eqalign{{\rm Re}(\bar{u}_2u_3)&=u_2^1u_3^1+
            u_2^2u_3^2+u_2^3u_3^3+u_2^4u_3^4=0,\cr
           {\rm Re}(\bar{u}_4u_5)&= u_4^1u_5^1+
u_4^2u_5^2+u_4^3u_5^3+u_4^4u_5^4=0.\cr}\leqno{\red.9}$$
But orthogonality also implies
$$|u_2|^2=|u_3|^2=|u_4|^2=|u_5|^2={1\over 4}.$$
Then we have
$$|{\rm Im}(\bar{u}_2u_3)|^2=|{\rm Im}(\bar{u}_4u_5)|^2={1\over 16},$$
and this contradicts the $U(1)$ moment map equation
$$2p_1{\rm Im}(\bar{u}_2u_3)+2p_2{\rm Im}(\bar{u}_4u_5)= 0$$
if $p_1\not=p_2$. Repeating the argument for the other choices
of $1<i<j<k<l$ gives the result under the hypothesis $p_1<p_2<p_3.$ 
Thus, at a vertex we must have $u_1\neq 0.$ This proves the ``neither"
part of the statement.

Now assume that $i=1<j<k<l$ and suppose either $j<k$ or $k<l$ is
a quaternionic pair. We can take $j=2, k=3$ and let $l$ be arbitrary.
Then the orthogonality of the corresponding $\bbb$-matrix
again forces 
$${\rm Re}(\bar{u}_2u_3)=0.$$
But again the $U(1)$ moment map constraint is
$${\rm Im}(\bar{u}_2u_3)=0$$
giving a contradiction. This proves the ``nor" part of the lemma.
It is now clear that the vertices must automatically satisfy the
$U(1)$ moment map constraint $\nu_\bfp(\bfu)=0;$ hence, they are all
diffeomorphic to $O(4).$ Moreover, a simple counting shows that there are
precisely eight vertices.  \hfill\za

Our analysis suggests the importance of the following strata:
$$\eqalign{S_0&=\{\bfu\in N_\nu(\bfp)~|~u_1=0\},\cr
           S_1&=\{\bfu\in N_\nu(\bfp)~|~u_1\neq 0\},\cr
           S_2&=\{\bfu\in N_\nu(\bfp)~|~\hbox{some quaternionic pair 
           $(u_{2i},u_{2i+1})$ vanishes}\},\cr
           S_3&=\{\bfu\in N_\nu(\bfp)~|~\hbox{no quaternionic pair 
           $(u_{2i},u_{2i+1})$ vanishes}\}.\cr} \leqno{\red.10}$$
Then Lemma \red.7 easily implies that

\noindent{\sc Corollary} \red.11: \tensl Let $0<p_1<p_2<p_3$. Then
\item{(i)} $S_0\cup S_1=S_2\cup S_3=N_\nu(\bfp).$
\item{(ii)} $S_0\cap S_1=S_2\cap S_3=\emptyset.$
\item{(iii)} $S_0\cap S_2=\emptyset.$
\item{(iv)} $S_2\subset S_1.$ 
\item{(v)} $S_0\subset S_3.$ \tenrm

Notice that (iii) fails if $p_i=p_j$ for some $i\neq j.$ In particular it
fails for the level set $N_\nu$ of \ks.12 in the previous section, and this is
the reason that the quotient $\calm$ of \ks.16 is not smooth. We now are ready
to give necessary conditions to guarantee a smooth quotient.

\noindent{\sc Lemma \red.12}: \tensl Let $\bfp=(p_1,p_2,p_3)\in(\bbz_+)^3$
be pairwise relatively prime.
Then the isotropy group of the $Sp(1)\times U(1)_{{\bf p}}$
action at every point of $S_1$ is the identity.
\tenrm
 
\noindent{\sc Proof}: The action of $Sp(1)\times U(1)_{{\bf p}}$ on $\bbh^7$
is the diagonal action of $Sp(1)$ by quaternionic multiplication by a unit
quaternion $\grl$ on the left, and the matrix multiplication $\bfu\mapsto
f_\bfp(t)\bfu$ for the $U(1)_\bfp$ action. These two actions clearly commute.  
Since $u_1\not=0$ we immediately get that $\lambda=1$. Consider
the set where a quaternionic pair $(u_{6},u_{7})=(0,0)$.
Then the fixed point equation becomes
$$A(p_1t)=A(p_2t)=\bbi_2,\leqno{\red.13}$$
which has only the trivial solution provided that 
${\rm gcd}(p_1,p_2)=1$.  Setting the other two quaternionic pairs
to be zero gives ${\rm gcd}(p_1,p_3)={\rm gcd}(p_2,p_3)=1$. As one cannot set
more than one quaternionic pair equal to $(0,0)$ the lemma is proved.
\hfill\za  

\noindent{\sc Lemma} \red.14: \tensl Let $\bfp=(p_1,p_2,p_3)\in(\bbz_+)^3$
satisfy the four conditions ${\rm gcd}(p_1\pm p_2,p_1\pm p_3)=1.$ Then the
isotropy group of the $Sp(1)\times U(1)_{{\bf p}}$  action at every point of
$S_0$ is the identity. \tenrm

\noindent{\sc Proof}: To determine the conditions for fixed points of the
action we consider the following equations
$$A(p_it)\pmatrix{u_{2i}\cr u_{2i+1}\cr}
=\pmatrix{a_i&b_i\cr-b_i&a_i\cr}\pmatrix{u_{2i}\cr u_{2i+1}\cr}=
\lambda \pmatrix{u_{2i}\cr u_{2i+1},\cr}\ \ \ \ i=1,2,3$$
for $\grl\in Sp(1)$ and $t\in [0,2\pi).$
For each $i=1,2,3$ this reads
$$a_i u_{2i}+b_i u_{2i+1}=\lambda u_{2i},$$
$$-b_i u_{2i}+a_i u_{2i+1}=\lambda u_{2i+1}.$$
For each $i=1,2,3$ we multiply the first equation from the right
by $\overline u_{2i}$ and the second by $\overline u_{2i+1}$ to get
$$a_i |u_{2i}|^2+b_i u_{2i+1}\overline u_{2i}=\lambda |u_{2i}|^2,$$
$$-b_i u_{2i}\overline u_{2i+1}+a_i |u_{2i+1}|^2=\lambda |u_{2i+1}|^2.
\leqno{\red.15}$$
By adding these two equations we get
$$a_i(|u_{2i}|^2+|u_{2i+1}|^2)+b_i(u_{2i+1}\overline u_{2i}-
u_{2i}\overline u_{2i+1})=\lambda(|u_{2i}|^2+|u_{2i+1}|^2),\ \ i=1,2,3.$$
By (iii) of Lemma \red.10 the term multiplying $\grl$ on the right hand side
of this equation never vanishes. This gives for each $i=1,2,3$ 
$$\eqalign{{\rm Re}(\lambda)&=a_i,\cr 
           {\rm Im}(\lambda)&= b_i{u_{2i+1}\overline u_{2i}-
u_{2i}\overline u_{2i+1}\over |u_{2i}|^2+|u_{2i+1}|^2}.\cr}
\leqno{\red.16}$$
The first of these equations gives
$$a_1=a_2=a_3,\leqno{\red.17}$$
and combining this with $a_i^2+b_i^2=1$ implies
$$b_1=\pm b_2=\pm b_3.\leqno{\red.18}$$
Let us write $\tau=e^{it}$. Then \red.17 and \red.18 give 
$$\tau^{p_1}=\tau^{\pm p_2},\qquad \tau^{p_1}=\tau^{\pm p_3}.\leqno{\red.19}$$ 
These have only trivial solutions if and only if   
${\rm gcd}(p_1\pm p_2,p_1\pm p_3)=1.$ 
\hfill\za

It is convenient to make the following:

\noindent{\sc Definition \red.20}: \tensl Let 
${\bf p}=(p_1,p_2,p_3)\in\bbz^3$.
We say that the weight vector ${\bf p}$ is {\it admissible} if
$0<p_1<p_2<p_3$, ${\rm gcd}(p_i,p_j)=1$ for all $i<j$,
and ${\rm gcd}(p_1\pm p_2,p_1\pm p_3)=1.$
\tenrm

It now follows immediately from Lemmas \red.12, \red.14 and Definition
\red.20 that

\noindent{\sc Theorem \red.21}: \tensl The $Sp(1)\times U(1)_{\bf p}$
action on $N_\nu(\bfp)$ is free if and only if ${\bf p}\in\bbz_+^3$ is 
admissible.
\tenrm

Note that there are infinitely many admissible weight vectors. For example
we can take ${\bf p}=(2k-1,2k,2k+1)$, where $k\in\bbz_+$. Thus, there
are infinite families of smooth quotients ${\cal M}({\bf p})$ and
infinite families of the associated triples ${\cal M}({\bf p})\rightarrow
{\cal Z}({\bf p})\rightarrow {\cal O}({\bf p})$ with their
(orbifold) Einstein metrics. 

Theorem A now follows from the Theorem
\red.21 and various theorems concerning 3-Sasakian (complex contact, 
quaternionic K\"ahler) reductions [BGM1, GL]. The last statement of Theorem A
follows from Corollary \ks.7.

We briefly return to the ${\bf p}=(1,1,1)$ case of
the previous section. As already observed
we get singularities here because the corresponding
action is not even locally free on the level set of the moment map.
But it is easy to see that it is quasi-free. That is there are
only two types of isotropy subgroups in the circle: the identity and the
whole group. In such a situation Dancer and Swann [DS] observed that 
3-Sasakian quotients, stratify as the union of
3-Sasakian manifolds [DS]. This is true in this case in particular, but as we
saw in the previous section the two strata nicely fit together
and one gets a compact 3-Sasakian orbifold. Let us explicitly
describe the singular part ${\cal M}_1\subset {\cal M}(1,1,1).$

Note that if $p_1=p_2=p_3=1$ then $a_1=a_2=a_3=a$
and $b_1=b_2=b_3=b$ and we can add equations (\red.15) to get
$$a+b\rho=\lambda,\leqno{\red.22}$$
where now
$$\rho=\sum_{i=1,2,3}
(u_{2i+1}\overline u_{2i}-
u_{2i}\overline u_{2i+1}).\leqno{\red.23}$$
When $u_1\not=0$ one does not get any fixed points of the action.
But when $u_1=0$ for any imaginary unit $\rho$
there is a $U(1)\subset Sp(1)\times U(1)$ subgroup
$$(\cos t+\rho \sin t, A(t))\in Sp(1)\times U(1),$$
which acts trivially on the following set
$$\{\bfu\in N_\nu ~|~ u_{2i+1}=\rho u_{2i}\}.
\leqno{\red.24}$$
In this case all 4 moment map equations reduce to the same one
and it reads
$$\sum_{i=1,2,3}\overline u_{2i}\rho u_{2i}=0.$$
For any fixed $\rho$ we can recognize this set as the
complex Stiefel manifold $U(3)/U(1)$ and it follows that the singular stratum
${\cal M}_1$ is precisely 
the quotient $SU(3)/U(1)={\cal S}(1,1,1).$

The geometry of the smooth families ${\cal M}({\bf p})$ is rather
interesting. First we observe that these spaces cannot be toric.
This can be seen in several different ways, for example by careful
analysis of the associated foliations. One can also generalize the
analysis of [BGM2] to show that the only isometries of the
level set of the moment map $N_\nu({\bf p})\subset S^{27}\subset\bbr^{28}$
can come from the restriction of the isometries of the Euclidean space 
$\bbr^{28}$. From this we conclude 

\noindent{\sc Theorem \red.25}: \tensl Let $\bfp$ be admissible so that
$\calm(\bfp)$ is a smooth compact 3-Sasakian 11-manifold. Then
the Lie algebra ${\rm Isom}^0({\cal M}(\bfp), g(\bfp))$
of the group of 3-Sasakian isometries of $\calm(\bfp)$
is isomorphic to
$\bbr^2\oplus {\Got s}{\Got p}(1)$.  In particular, all such
quotients are non-toric.
\tenrm

This proves part of Theorem C of the introduction which relates to the
11-dimensional quotients. 

Finally, observe that ${\cal M}({\bf p})$
contains a special 7-manifold  which is embedded as a 3-Sasakian
submanifold. Define
$$S_0({\bf p})={N_\nu({\bf p})\cap \{u_1=0\}\over Sp(1)\times U(1)_{\bf p}}.
\leqno{\red.26}$$
One can see that $S_0({\bf p})$ is a submanifold and, as it is
itself a 3-Sasakian reduction, it must be a 3-Sasakian submanifold.
One can even identify this space. Observe that the  classical
group isomorphism ${\rm Spin}(6)\simeq SU(4)$ implies that the
$Sp(1)$ quotient yields the homogeneous 3-Sasakian 11-dimensional manifold
${\cal S}(1,1,1,1)$. Hence, $S_0({\bf p})$ is either a $U(1)_{\bf p}$-reduction
of ${\cal S}(1,1,1,1)$ or, equivalently, a $T^2$-reduction of
$S^{15}$. Hence, there exists an admissible integer weight matrix 
$\Omega\in{\cal A}_{2\times 4}(\bbz)$ (see [BGMR]) such that
$S_0({\bf p})\simeq {\cal S}(\Omega)$.  Hence $S_0({\bf p})$ is toric with
second Betti number equal to 2.

\bigskip
\def\cals{{\cal S}}

\centerline {\bf \low. Further Reductions of ${\cal M}(p_1,p_2,p_3)$
by a Circle} 
\medskip
In this section we will examine reductions of ${\cal M}(p_1,p_2,p_3)$
by isometric circle actions. More generally we shall consider an 
arbitrary 2-torus subgroup of the maximal 
torus $T^2\subset T^3\subset SO(7)$.
Let
$${\Theta}=\pmatrix{p_1&p_2&p_3\cr
                    q_1&q_2&q_3\cr}\in {\cal M}_{2\times3}(\bbz)
\leqno{\low.1}$$
be an arbitrary integral $2\times3$ matrix and
define the homomorphism $f_\Theta:T^2\rightarrow SO(7)$
$$f_{\Theta}(t,s)=\pmatrix{1&0&0&0\cr 0&A(p_1t+q_1s)&0&0\cr 0&0&A(p_2t
+q_2s)&0\cr
0&0&0&A(p_3t+q_3s)\cr}\in SO(7),\leqno{\low.2}$$
where $A$ is the $SO(2)$ rotations as in \red.3. Note that
if $\bfp=(p_1,p_2,p_3)$ were admissible then we would be considering
arbitrary isometric circle actions on the quotient ${\cal M}(\bfp)$;
however, so far we assume nothing about $\Theta$.
The homomorphism $f_{\Theta}(t,s)$
yields a 2-torus action on the homogeneous 3-Sasakian manifold
$SO(7)/SO(3)\times Sp(1)$ via left multiplication
$f_{\Theta}(t,s){\bf u}$ and the moment map can now be written as
$$\nu_{\Theta}(u_1,...,u_7)=\pmatrix
{\sum_{\alpha=1,2,3}p_\alpha(
\overline u_{2\alpha}u_{2\alpha+1}-\overline u_{2\alpha+1}u_{2\alpha})\cr
\sum_{\alpha=1,2,3}q_\alpha(
\overline u_{2\alpha}u_{2\alpha+1}-\overline u_{2\alpha+1}u_{2\alpha})\cr}
\in\bbr^2\otimes{\Got s}{\Got p}(1).
\leqno{\low.3}$$

\noindent
We begin our analysis of this quotient by considering the level set 
of the moment map. 

\noindent{\sc Definition \low.4}: \tensl Let $N_\nu(\Theta)\subset S^{27}$
denote the level set of the 3-Sasakian moment map of the 
$Sp(1)\times T^2_\Theta$-action, {\it i.e.,} 
$N_\nu({\Theta})\equiv N\cap \{\nu_{{\Theta}}^{-1}({\bf 0})\}$ and let
$${\cal M}(\Theta)=
N_\nu({\Theta})/
Sp(1)\times T^2_{\Theta}. \leqno{\low.5}$$
\tenrm
We want to determine for which $\Theta\in{\cal M}_{2\times3}(\bbz)$
the 7-dimensional quotient ${\cal M}(\Theta)$ is an orbifold and,
if possible 
which weight matrices yield smooth quotients. Let us define
$$\Delta_{ij}=
\Delta_{ij}(\Theta)={\rm det}\pmatrix{p_i&p_j\cr q_i& q_j\cr}, \qquad
1\leq i<j\leq3, \leqno{\low.6}$$ 
the three minor determinants of $\Theta$.

\noindent{\sc Lemma} \low.7: \tensl The action of
$Sp(1)\times T^2_{\Theta}$ on $N_\nu({\Theta})\cap\{u_1\not=0\}$ is
\item{(i)} locally free if and only if 
$\Delta_{ij}(\Theta)\not=0,\qquad\forall~ 1\leq i<j\leq3$,
\item{(ii)} free if and only if $|\Delta_{ij}(\Theta)|=1,\qquad\forall~ 
1\leq i<j\leq3$.
\tenrm

\noindent{\sc Proof}: Since $u_1\not=0$ we must have $\lambda=1$ and, hence,
it is enough to consider the $T^2_\Theta$-action. Now, suppose
that a quaternionic pair, say $(u_6,u_7)$ vanishes. Then the
fixed point equation reads:
 
$$A(p_it+q_is)\pmatrix{u_{2i}\cr u_{2i+1}\cr}
=\pmatrix{u_{2i}\cr u_{2i+1}\cr},\ \ \ \ i=1,2\leqno{\low.8}$$
for $t,s\in [0,2\pi),$ or, equivalently,
$$A(p_it+q_is)=\bbi_2,
\ \ \ \ i=1,2.\leqno{\low.9}$$
Let $\tau=e^{it}$ and $\rho=e^{is}$. Then we can rewrite \low.9
as
$$\tau^{p_i}\rho^{q_i}=1, \ \ \ \ i=1,2.\leqno{\low.10}$$
This has only discrete solutions provided $\Delta_{12}(\Theta)\not=0$.
Furthermore, the isotropy group at all such points will be trivial
provided $\Delta_{12}(\Theta)=\pm1$. This proves the lemma.
\hfill\za

Note that the second condition is already very restrictive as it
says that any $2\times2$ submatrix of $\Theta$
must be an element of $PSL(2,\bbz)$. However, there are many matrices which
satisfy both conditions, for example
$$\Theta_1=\pmatrix{1&0&1\cr0&1&1}, \ \ \ \ \ 
\Theta_2=\pmatrix{9&2&7\cr 40&9&31\cr}.$$
It remains to analyze the fixed point equations on
$N_\nu({\Theta})\cap\{u_1=0\}$. We now prove 

\noindent{\sc Lemma} \low.11: \tensl The action of
$Sp(1)\times T^2_{\Theta}$ on $N_\nu({\Theta})$ is
locally free if and only if 
$\Delta_{ij}(\Theta)\not=0,\ \ \forall~ 1\leq i<j\leq3$ and
$${\square_\mp^\mp}={\rm det}\pmatrix{p_1\mp p_2&q_1\mp q_2\cr
           p_1\mp p_3&q_1\mp q_3\cr}\not=0 \leqno{\low.12}$$
Furthermore, there is no weight matrix $\Theta\in{\cal M}_{2\times3}(\bbz)$ 
for which the actions is free. 
\tenrm

\noindent{\sc Proof}: First let us clarify that what we mean in
\low.12 is that four determinants $\square_+^+$,
$\square_+^-$, $\square_-^+$, and $\square_-^-$ must vanish 
(in any row we can choose either upper or lower signs). Since the
$\Delta_{ij}(\Theta)\not=0$, by the previous
Lemma, we know that the action is locally free on the
$u_1\not=0$ part. Hence, it is enough to consider $u_1=0$. Here the
analysis is similar to the one presented in Lemma \red.14 and 
it is entirely based on the fact that no 
quaternionic pair can vanish. As a result
we get the following analogue of the fixed point equations
\red.19:
$$\tau^{p_1}\rho^{q_1}=\bigl(\tau^{p_2}\rho^{q_2}\bigr)^{\pm1},
\qquad \tau^{p_1}\rho^{q_1}=\bigl(\tau^{p_3}\rho^{q_3}\bigr)^{\pm1}.
\leqno{\low.13}$$ 
We can rewrite these as
$$\tau^{p_1\mp p_2}=\rho^{-q_1\pm q_2},\qquad
\tau^{p_1\mp p_3}=\rho^{-q_1\pm q_3}. \leqno{\low.14}$$
These are four systems of two equations in $(\tau,\rho)$ variables. We want 
all four of them to have at most discrete solutions. This requires that the
four determinants
$${\rm det}\pmatrix{p_1\mp p_2&-q_1\pm q_2\cr
           p_1\mp p_3&-q_1\pm q_3\cr}=-{\rm det}
\pmatrix{p_1\mp p_2&q_1\mp q_2\cr
           p_1\mp p_3&q_1\mp q_3\cr}\leqno{\low.15}$$
do not vanish and gives \low.12.

The fact that orbifold singularities are always present requires more
subtle analysis. To have smooth quotients we must assume
$$\Delta_{ij}(\Theta)=\pm1, \qquad \forall~ 1\leq i<j\leq3,
\leqno{\low.16}$$
one one hand, and
$$(\square_\mp^\mp)=\left|\matrix{p_1\mp p_2&q_1\mp q_2\cr
           p_1\mp p_3&q_1\mp q_3\cr}\right|=\pm1 \leqno{\low.17}$$
on the other. A simple computation relates all of
these four determinants to the 3 minor
determinants $\Delta_{ij}(\Theta)$ and we get
$$\square_-^-={\rm det}\pmatrix{p_1- p_2&q_1- q_2\cr
           p_1- p_3&q_1- q_3\cr}=\Delta_{12}(\Theta)+
\Delta_{23}(\Theta)-\Delta_{13}(\Theta),$$
$$\square_+^-={\rm det}\pmatrix{p_1- p_2&q_1- q_2\cr
           p_1+ p_3&q_1+ q_3\cr}=\Delta_{12}(\Theta)-
\Delta_{23}(\Theta)+\Delta_{13}(\Theta),$$
$$\square_-^+={\rm det}\pmatrix{p_1+ p_2&q_1+ q_2\cr
           p_1- p_3&q_1- q_3\cr}=-\Delta_{12}(\Theta)-
\Delta_{23}(\Theta)-\Delta_{13}(\Theta),$$
$$\square_+^+={\rm det}\pmatrix{p_1+ p_2&q_1+ q_2\cr
           p_1+ p_3&q_1+ q_3\cr}=-\Delta_{12}(\Theta)+
\Delta_{23}(\Theta)+\Delta_{13}(\Theta).\leqno{\low.18}$$
Now, because of \low.16, all three minor determinants must be $\pm1$.
This gives 8 possible combinations of the values of $\Delta_{ij}(\Theta)$.
It is trivial to check that for any one out of these eight at least
two determinants $\square_\mp^\mp$ will be equal to $\pm3$ (the other
six all being equal to $\pm1$). Hence, even if we choose $\Theta$ so
that \low.16 holds, the quotient will necessarily have
orbifold singularities of type $\bbz_3$. This concludes the proof 
of Lemma \low.11\hfill\za 

Using the calculation in the proof of the above lemma we restate  
condition \low.12 to get

\noindent{\sc Theorem} \low.19: \tensl The action of
$Sp(1)\times T^2_{\Theta}$ on $N_\nu({\Theta})$ is
locally free if and only if 
\item{(1)} all their determinants
$\Delta_{12}(\Theta), \Delta_{23}(\Theta), \Delta_{13}(\Theta)$ do not
vanish, and
\item{(2)} their sum does not vanish, and
\item{(3)} none of the determinants is equal to the sum of the other two.

\noindent
In such a case the quotient ${\cal M}(\Theta)$ is a compact 7-dimensional
3-Sasakian orbifold. Furthermore, there is no weight matrix
$\Theta$ for which ${\cal M}(\Theta)$ is a smooth manifold.
\tenrm

Now, Theorem C of the introduction follows from Theorem \low.19,
and the fact that the fundamental 
3-dimensional foliation ${\cal M}(\Theta)\rightarrow{\cal O}(\Theta)$, 
in the case ${\cal M}(\Theta)$ is a compact orbifold,
yields a compact self-dual Einstein orbifold with a
positive scalar curvature (orbifold) metric as the space of leaves. 
The fact that this
metric has only one Killing vector field follows from an appropriate
generalization of Theorem B.

\noindent{\sc Remark} \low.20: Note that both ${\cal M}(\Theta)$
and ${\cal O}(\Theta)$ depend only on the three minor determinants
$\Delta_{12}(\Theta), \Delta_{23}(\Theta), \Delta_{13}(\Theta)$ rather
that on $\Theta$ itself. Different weight matrices can certainly lead
to equivalent quotients. One could compute the self-dual
Einstein metrics $g(\Theta)$ on ${\cal O}(\Theta)$ explicitly. Locally we can
change variables so that
$$\pmatrix{1&0&0&0\cr 0&A(p_1t\!+\!q_1s)&0&0\cr 0&0&A(p_2t\!+q_2s)&0\cr
0&0&0&A(p_3t\!+\!q_3s)\cr}\!=\!
\pmatrix{1&0&0&0\cr 0&A(\lambda_1)&0&0\cr 0&0&A(
\lambda_2)&0\cr
0&0&0&A(a\lambda_1\!+\!b\lambda_2)\cr},
$$
where $a=-{\Delta_{23}(\Theta) \over \Delta_{12}(\Theta)}$ and
$b={\Delta_{13}(\Theta) \over \Delta_{12}(\Theta)}$ are now non-zero rational
numbers. Such a choice simplifies both the 2-torus action as well
as the moment map equations. The self-dual Einstein 
quotient metric in question (up to scale) will depend on these
two parameters $g(\Theta)=g(a,b)$.

\bigskip 
\def\tnu{\tilde{\nu}} 
 
\centerline {\bf \hig. 15-Dimensional Examples} 
\medskip 
 
At first, it may appear that there are obvious higher dimensional  
analogues of our construction. However, a simple parity argument 
shows that such actions do not yield 3-Sasakian metrics other than 
in dimensions 7, 11, and 15. One can have orbifold metrics in any  
3-Sasakian dimension. On the other hand the existence of 
these smooth quotient in dimension 11 and 15 is closely related to 
the geometry of $G_2$ and ${\rm Spin}(7)$. 
 
We consider two separate cases of $U(1)_{\bf p}\subset SO(2k+1)$  
and $U(1)_{\bf p} \subset SO(2k)$. 
In the $SO(2k+1)$ case it suffices to take $k=4$. Then 
${\bf p}=(p_1,p_2,p_3,p_4)$ and the new $U(1)_{\bf p}$ action is 
defined by adding an $SO(2)$ matrix $A(p_4t)$ and let it act by 
rotating the additional quaternionic coordinates $(u_8,u_9)$. To get 
a smooth quotient any triple 
$p_i<p_j<p_k$ would have to be admissible according to the 
definition \red.31, as is easily seen from the analysis given in Section 2. 
However, this is clearly impossible as admissibility implies 
that each triple contains two odd and one even number. 
 
In the $SO(2k)$ case we have already seen in the previous section that 
$k=3$ leads to the toric 7-manifolds ${\cal S}(\Omega)$ with $b_2=2$. 
Now we will show that, in fact, $k=4$ leads to  
smooth quotients.  
Let 
${\bf p}=(p_1,p_2,p_3,p_4)\in \bbz^4$ and define the following homomorphism 
$$f_{{\bf p}}(t)=\pmatrix{A(p_1t)&0&0&0\cr 0&A(p_2t)&0&0\cr 0&0&A(p_3t)&0\cr 
&0&0&A(p_4t)\cr}\in SO(8),\leqno{\hig.1}$$ 
where 
$$A(p_it)=\pmatrix{\cos(p_it)&\sin(p_it)\cr-\sin(p_it)&\cos(p_it)\cr} 
\in SO(2),\ \ \ \ \ i=1,2,3,4 \leqno{\hig.2}$$ 
are 2-dimensional real rotations. As before we can choose 
all weights to be non-negative. Further note that at most 
one of the weights can vanish. We prove the following  
 
\noindent{\sc Lemma \hig.3}: \tensl Let  
$(p_1,p_2,p_3,p_4)\in\bbz^4$. Then the action of $U(1)_{{\bf p}}\times Sp(1)$ 
on $N_\nu(\bfp)$ is free if and only if  
$0\leq p_1<p_2<p_3<p_4$,  
${\rm gcd}(p_i,p_j,p_k)=1$ and  
${\rm gcd}(p_i\pm p_j, p_i\pm p_k)=1$ for any triple in $\bfp$. 
\tenrm 
 
\noindent{\sc Proof}: First, let us assume that all of the weights are 
non-negative. Consider a triple, say $(p_1,p_2,p_3)$. Set $u_7=u_8=0.$ 
Then the analysis of the previous section shows that the three 
must be distinct and that we must have  
${\rm gcd}(p_1\pm p_2, p_1\pm p_2)=1$. However, we no longer need 
the three weights to be pairwise relatively prime as one cannot 
set two of the quaternionic pairs $(u_{2i-1}, u_{2i}), i=1,2,3,4$ equal to 
$(0,0)$ at the same time. One such quaternionic pair can vanish; hence,  
we need ${\rm gcd}(p_1,p_2,p_3)=1$ to get a free action. 
The analysis in the case when $p_1=0$ is similar. 
Then one sees that the triple $(p_2,p_3,p_4)$ 
has to be admissible in the sense of the Definition 
\red.31. But that is what Lemma \red.9 says in this case. 
\hfill\za 
 
It is immediately clear that one cannot extend this construction 
for $k>4$ without admitting orbifold singularities in the quotient. 
 
Theorem B follows from the above lemma except for its last statement. 
When $\bfp=(1,1,1,1)$ the quotient is singular. It easy to see 
that the action is quasi-regular just as in the 11-dimensional case.  
The fact that the two strata fit together giving the 3-Sasakian orbifold 
$\bbz_2\backslash {\rm Spin}(7)/{\rm Spin}(4)$ can be seen by first 
identifying the zero level  
set $N_\nu(1,1,1,1)$ with the $U(1)\cdot \big({\rm Spin}(7)/Sp(1)\big)$  
and observing 
that $U(1)\cap {\rm Spin}(7)\simeq U(1)\cap SU(4)\simeq\bbz_2$ [BGOP]. 
 
Actually, $\bfp=(1,1,1,1)$ is not quite the reduction considered 
by Kobak and Swann [KS3]. Instead they use $\bfp=(0,0,0,1)$ for the 
$U(1)$-action. The latter is easily seen to give ${\rm Gr}_4(\bbr^7)/\bbz_2$ 
as the quaternionic K\"ahler quotient. This is a simple consequence 
of the fact that a "zero momentum" hyperk\"ahler reduction of $\bbh^2$ by any 
circle action is isometric to $\bbh/\Gamma$ ($\Gamma=\bbz_2$ for the 
standard action). The identification of the cases $\bfp=(1,1,1,1)$ 
and $\bfp=(0,0,0,1)$ owes to the isomorphism between 
${\Got s}{\Got o}(6)$ and ${\Got s}{\Got u}(4)$. 
 
The fact that these quotients cannot be toric follows from the same type 
of argument as the one used for 11-dimensional quotients. The Lie algebra of 
the of the group of the 3-Sasakian isometries is $\bbr^3\oplus  
{\Got s}{\Got p}(1)$. 
 
The smooth 15-dimensional manifolds 
${\cal M}(\bfp)$ for $p_1=0$ and $p_1>0$ are geometrically different. 
In the first case ${\cal M}(\bfp)$ contains two copies of the 11-dimensional 
3-Sasakian manifold ${\cal M}(p_2,p_3,p_4)$ which intersect in the 
7-dimensional toric 3-Sasakian submanifold $\cals(\Omega(p_2,p_3,p_4))$. 
When $p_1>0$ then ${\cal M}(\bfp)$ does not have any obvious 11-dimensional 
3-Sasakian submanifolds. However, we do get 4 disjoint toric 7-dimensional 
3-Sasakian submanifolds $\cals(\Omega(p_1,\ldots,\hat p_i,\ldots,p_4))$ 
by setting one of the quaternionic pairs $(u_{2i-1},u_{2i})=(0,0)$. 
 
We can also consider non-zero momentum  
$\xi\in {\Got s}{\Got p}(1)$-deformations of 
the hyperk\"ahler metric on the cones $C({\cal M}({\bf p}))$. (Up to scale 
one can set $\xi$ to be any imaginary quaternion, so we really have just 
one parameter family). 
In some sense they are all deformations of Kronheimer's 
hyperk\"ahler metrics on the two nilpotent orbits in 
${\Got s}{\Got l}(3,\bbc)$ and ${\Got s}{\Got o}(6, \bbc)$. 
Unfortunately such hyperk\"ahler quotients are rarely free from 
orbifold singularities. In fact, in the 12-dimensional case we never 
get any complete metrics. In 16 dimensions we do get a complete metric 
only when ${\bf p}=(1,1,1,1)$ (or equivalently ${\bf p}=(0,0,0,1)$). 
The metric is $SU(4)$-invariant and it gives the 
Kronheimer metric on the 16-dimensional nilpotent orbit 
of ${\Got s}{\Got o}(6, \bbc)$ in the $\xi\rightarrow0$ scaling limit. 
More generally, cohomogeneity 2 metrics were studied by Kobak and Swann 
[KS4]. 
 
As each classical nilpotent orbit is a hyperk\"ahler reduction of 
some quaternionic vector space [KS2] it is tempting to undertake a more 
systematic study of the following problem: Which of the nilpotent 
orbits can give rise to compact 3-Sasakian manifolds? Certainly, any time a 
quotient involves some $U(1)$-factor one can introduce weights. 
However, as demonstrated here, requiring smoothness often puts very 
severe restrictions on the weights. We plan to address some of these 
questions in a future work [BGOP]. 
 
\noindent{\sc Remark \hig.4}: Note that all of the quotients considered in this 
paper are examples of toric reductions of the 3-Sasakian homogeneous space 
associated to the real Grassmanian $Gr_4(\bbr^n)$. This space 
is $SO(n)/SO(n-4)\times Sp(1)$ and it has $SO(n)$ as the group 
of isometries preserving the 3-Sasakian structure. Consider 
the maximal torus $T^l\subset SO(n)$. Then the relevant question is: 
Which subgroups $T^m\subset T^l$ yield smooth 3-Sasakian quotients? 
Sections \red\  and \low\  give the complete analysis of the $n=2k+1=7$ 
case. In Section \red\ $m=1$ and in Section \low\  $m=2$ which exhaust 
all interesting possibilities. This section considers $n=2k=8$ case 
with $m=1$. One can see that $m>1$ does not yield any smooth quotient 
but we leave the analysis of this to a future work where we plan to 
give a complete answer in the most general case of arbitrary $(m,n)$ 
[BGP]. 

\bigskip
\def\tnu{\tilde{\nu}}

\centerline {\bf \top. Comments on the Topology of $\calm(\bfp)$ and Related
Spaces} \medskip

In this section we denote by $\calm(\bfp)$ either one of the 11 or 15
dimensional 3-Sasakian manifolds discussed in sections \red~ and \hig~ with
$\bfp$ admissible. Actually as discussed below $\calm(\bfp)$ will denote a
component of the manifolds discussed previously. It would be interesting to
know the topology of our quotients, most importantly $\pi_1(\calm(\bfp))$ and
$H_2({\cal M}({\bf p}),\bbz)$.  For this one needs to understand the topology
of the level set of the moment map $N_\nu(\bfp).$ Of course, we do know that
$\pi_1(\calm(\bfp))$ is finite and that the odd Betti numbers of $\calm(\bfp)$
vanish up to the middle dimension [GS]. However, beyond this not much explicit
topological information can be obtained. For example, so far we have been
unable to determine whether $\calm(\bfp)$ and $N_\nu(\bfp)$ are even
connected. This presents no real problem as we shall always mean by these
spaces connected components such that $N_\nu(\bfp)$ is a $S^1\times Sp(1)$
bundle over $\calm(\bfp).$ Generally, the determination of the topology of an
intersection of real quadrics such as $N_\nu(\bfp)$ is quite complicated. The
analysis in previous work [BGM1] and [BGMR] relied heavily on very specialized
information. In the former case the level sets in question were diffeomorphic
to certain Stiefel manifolds whose topology is completely understood, whereas
in the later case there was a large symmetry group whose quotient was a two
dimensional space whose topology could be analysized. In the present case we
meet with no such good fortune. Nevertheless, we are able to obtain a small
amount of general information about the topology of $N_\nu(\bfp).$ Our first
result is that for $\bfp$ and $\bfp'$ admissible the level sets $N_\nu(\bfp)$
and $N_\nu(\bfp')$ are diffeomorphic. This does not hold generally as the
case of $N_\nu(1,1,1)$ shows. In this case the Jacobian matrix drops rank
making the quotient singular.  

\noindent{\sc Lemma} \top.1: \tensl For $\bfp$ and $\bfp'$ admissible, 
the level sets $N_\nu(\bfp)$ and $N_\nu(\bfp')$ are diffeomorphic. \tenrm 
 
\noindent{\sc Proof}: For $\bfp$ admissible there are no fixed points of the 
$U(2)$ action so by the 3-Sasakian version of a well known result in 
symplectic geometry, zero is a regular value of the $U(2)$ moment map  
$$\tilde{\nu}_\bfp: S^{27}\ra{1.3} \bbr^{12}=\bbr^3\times\bbr^3\times  
\bbr^3\times\bbr^3$$  
defined by $\tnu_\bfp =(\nu_\bfp,\mu_i,\mu_j,\mu_k).$ Now the level sets 
$N_\nu(\bfp)$ are defined for all $\bfp\in \bbr^3.$  Moreover, for any 
non-zero $\rho \in \bbr$ we see that level sets of $\tnu_{\rho\bfp}$ and 
$\tnu_\bfp$ coincide, that is $N_\nu(\rho\bfp)=N_\nu(\bfp).$   Thus, by scaling 
we can choose $\rho$ such that $\tnu_{\rho\bfp'}$ is in an $\gre$-neighborhood 
of $\tnu_{\rho\bfp}$ in $C^\infty(S^{27},\bbr^{12})$ with the $C^\infty$ 
compact-open topology. Since zero is a regular value of both  
$\tnu_{\rho\bfp}$ and $\tnu_{\rho\bfp'},$ it follows from a well known theorem 
(cf. [BCR], 14.1.1) that $N_\nu(\bfp)$ and $N_\nu(\bfp')$ are diffeotopic, that 
is there is a one parameter family of diffeomorphisms $\phi_t:S^{27}\ra{1.4} 
S^{27}$ parameterized by the unit interval such that $\phi_0$ is the identity 
and $\phi_1$ takes $N_\nu(\bfp')$ diffeomorphically to $N_\nu(\bfp).$ \hfill\za 

This immediately implies that the homotopy groups as well as the cohomology
rings of $N_\nu(\bfp)$ and $N_\nu(\bfp')$ are isomorphic. 

Now once and for all we shall choose a component of $N_\nu(\bfp)$ and the
corresponding component of $\calm(\bfp)$ so that $N_\nu(\bfp)$ is an
$S^1\times Sp(1)$ bundle over $\calm(\bfp)$ with both base space and total
space connected. We also denote by $\calh(\bfp)$ the circle bundle over
$\calm(\bfp)$ that coincides with the quotient of the corresponding component
of $N_\nu(\bfp)$ by the $Sp(1)$ action. We now study the latter as a 
fibration, namely $S^3\ra{1.3} N_\nu(\bfp)\fract{\pi}{\ra{1.3}} \calh(\bfp).$
Note that since $\phi_1$ in Lemma \top.1 is not necessarily a bundle map, we
cannot claim that the manifolds $\calh(\bfp)$ and $\calh(\bfp')$ are
diffeomorphic. Nevertheless, we can obtain some useful information about their
cohomology and homotopy groups. The manifolds $\calh(\bfp)$ are of interest in
their own right since as discussed briefly in the next section they admit
hypercomplex structures.

Consider the following commutative diagram of Gysin sequences with $\bbz$
coefficients:
$$\matrix{H^{r+3}(N_\nu(\bfp))\ra{1.1}&H^r(\calh(\bfp))\fract{\cup
\chi}{\ra{1.1}}&H^{r+4}(\calh(\bfp))\fract{\pi^*}{\ra{1.1}}
&H^{r+4}(N_\nu(\bfp))\ra{1.1} &H^{r+1}(\calh(\bfp))\cr
\decdnar{}&{}&&\decdnar{}& \cr
H^{r+3}(N_\nu(\bfp'))\ra{1.1}&H^r(\calh(\bfp'))\fract{\cup
\chi}{\ra{1.1}}&H^{r+4}(\calh(\bfp'))\fract{\pi^*}{\ra{1.1}}&H^{r+4}(N_\nu(\bfp'))      
\ra{1.1}&H^{r+1}(\calh(\bfp')),}\leqno{\top.2}$$
where the two vertical arrows are isomorphisms by Lemma \top.1, and $\cup
\chi$ denotes cupping by the Euler class of the bundle. The idea is to
construct, as best as possible, the missing vertical arrows and
relate the cohomology of $\calh(\bfp)$ and $\calh(\bfp')$ by the Five Lemma.
First we notice that setting $r=-3$ and $-2$ and using Lemma \top.1 gives
isomorphisms $$H^1(\calh(\bfp),\bbz)\approx H^1(N_\nu(\bfp),\bbz)\approx
H^1(N_\nu(\bfp'),\bbz)\approx H^1(\calh(\bfp'),\bbz)$$
$$H^2(\calh(\bfp),\bbz)\approx H^2(N_\nu(\bfp),\bbz)\approx
H^2(N_\nu(\bfp'),\bbz)\approx H^2(\calh(\bfp'),\bbz)$$
Next by setting $r=-1$ we have
$$\matrix{0&\ra{1.5}&H^{3}(\calh(\bfp))\fract{\pi^*}{\ra{1.5}}
&H^{3}(N_\nu(\bfp))&\ra{1.5}&H^0(\calh(\bfp))\cr
          \decdnar{=}&&&\decdnar{\phi^*}&&\decdnar{\approx} \cr
          0&\ra{1.5}&H^{3}(\calh(\bfp'))\fract{(\pi')^*}{\ra{1.5}}
&H^{3}(N_\nu(\bfp'))&\ra{1.5}&H^0(\calh(\bfp'))}\leqno{\top.3}$$
Since both $\pi^*$ and $(\pi')^*$ are injective, $\phi^*$ is an isomorphism,
and the diagram is exact and commutative,  we can define the missing vertical
arrow by $\psi =((\pi')^*)^{-1}\circ \phi^*\circ \pi^*$ and it is an
isomorphism. We thus have
$$H^{3}(\calh(\bfp),\bbz)\approx H^{3}(\calh(\bfp'),\bbz).$$
Now generally we cannot construct the missing vertical maps; however, we can
construct them if the groups are free. We thus change to rational coefficients
$\bbq.$

Consider now diagram \top.2 for $r=0.$ We can now fill in the second and last
columns with vertical arrows that are isomorphisms. Now considering the
diagram  with rational coefficients, we can split the middle groups, for all
admissable $\bfp,$ as
$$H^4(\calh(\bfp),\bbq)\approx \hbox{im}(\cup\chi) \oplus
\hbox{coker}(\cup\chi).$$ 
Choosing bases for these groups we can define the middle vertical map simply
by sending a basis element in
$\hbox{im}(\cup\chi)\subset H^4(\calh(\bfp),\bbq)$ to a basis element in
$\hbox{im}(\cup\chi)\subset H^4(\calh(\bfp'),\bbq),$ and a basis element in
$\hbox{coker}(\cup\chi)\subset H^4(\calh(\bfp),\bbq)$ to a basis element in
$\hbox{coker}(\cup\chi)\subset H^4(\calh(\bfp'),\bbq).$ That these subspaces
have the same dimension making this possible follows from exactness and
commutativity of the diagram. It then follows from the Five Lemma that this
middle arrow is an isomorphism. This argument is general and using a simple
induction we arrive at

\noindent{\sc Theorem} \top.4: \tensl For admissible $\bfp$ and $\bfp'$ we have
\item{(i)} $H^r(\calh(\bfp),\bbz)\approx H^r(\calh(\bfp'),\bbz)$ for
$r=0,1,2,3.$
\item{(ii)} $b_r(\calh(\bfp))=b_r(\calh(\bfp'))$ for all $r.$

\noindent In particular, $\calh(\bfp)$ and $\calh(\bfp')$ have isomorphic
rational cohomology groups. \tenrm

Here $b_r(M)$ denotes the rth Betti number of $M.$ Similar arguments can be
used for the homotopy groups, and combining these results with
well known facts about circle bundles over 3-Sasakian manifolds [BGM3] we
obtain

\noindent{\sc Proposition} \top.5: \tensl  For $\bfp$ and $\bfp'$ admissible,
we have isomorphisms: 
\item{(i)} $\pi_i(\calh(\bfp))\approx \pi_i(N_\nu(\bfp))\approx
\pi_i(\calh(\bfp'))$ for $i=0,1,2.$
\item{(ii)} $H^1(\calh(\bfp),\bbz)\approx H^1(N_\nu(\bfp),\bbz)\approx 
H^1(\calh(\bfp'),\bbz)\approx 0~\hbox{or}~ \bbz.$
\item{(iii)} $H^2(\calh(\bfp),\bbz)\approx H^2(N_\nu(\bfp),\bbz).$ \tenrm

Next we consider our manifolds of primary interest, namely $\calm(\bfp).$
First, it is easy to see the following relations:

\noindent{\sc Proposition} \top.6: \tensl For admissible $\bfp$ we have:
\item{(i)} $\pi_i(\calh(\bfp))\approx \pi_i(\calm(\bfp))$ for all $i>2.$
\item{(ii)} $b_2(N_\nu(\bfp))=\displaystyle{\cases{b_2(\calm(\bfp))-1 &if 
                               $b_1(N_\nu(\bfp))=0;$\cr                      
                         b_2(\calm(\bfp)) &if $b_1(N_\nu(\bfp))=1.$ \cr}}$ 
\tenrm

Next we have the analogue of Theorem \top.4 for our $3$-Sasakian manifolds
$\calm(\bfp).$

\noindent{\sc Theorem} \top.7: \tensl For admissible $\bfp$ and $\bfp'$ we have
\item{(i)} $H^r(\calm(\bfp),\bbz)\approx H^r(\calm(\bfp'),\bbz)$ for
$r=0,1,2,3.$
\item{(ii)} $b_r(\calm(\bfp))=b_r(\calm(\bfp'))$ for all $r.$

\noindent In particular, $\calm(\bfp)$ and $\calm(\bfp')$ have isomorphic
rational cohomology groups. \tenrm

\noindent{\sc Proof}: The proof of this Theorem is analagous to the proof
of Theorem \top.4 with diagram \top.2 replaced by the following  commutative
diagram of Gysin sequences:
$$\matrix{H^{r+1}(\calh(\bfp))\ra{1.1}&H^r(\calm(\bfp))\fract{\cup
\chi}{\ra{1.1}}&H^{r+2}(\calm(\bfp))\ra{1.1}&H^{r+2}(\calh(\bfp))\ra{1.1}
&H^{r+1}(\calm(\bfp))\cr
\decdnar{}&&&\decdnar{}& \cr
H^{r+1}(\calh(\bfp'))\ra{1.1}&H^r(\calm(\bfp'))\fract{\cup
\chi}{\ra{1.1}}&H^{r+2}(\calm(\bfp'))\ra{1.1}&H^{r+2}(\calh(\bfp'))\ra{1.1}
&H^{r+1}(\calm(\bfp')),}$$ 
where we have used Theorem \top.4 and its proof to construct the 
isomorphisms indicated by the vertical arrows.
\hfill\za

\noindent{\sc Remark} \top.8: Actually we can weaken the hypothesis that
$\bfp$ be admissible by noting any of the results of this section concerning
rational cohomology hold in the cases when $\calm(\bfp)$ is an orbifold
obtained as the quotient by a locally free action. It follows from the
analysis in section \red, that the action is locally free precisely when the
components of $\bfp$ are all distinct, and in this case the level sets
$\calh(\bfp)$ are smooth manifolds; hence, Theorem \top.4 and Proposition
\top.5 hold in this case as well. It is interesting to note that in the case
that $\calm(\bfp)$ is an orbifold, but not a smooth manifold, the smooth
manifold $\calh(\bfp)$ cannot be the trivial V-bundle. The above remarks
apply equally as well to the 7-dimensional orbifolds $\calm(\Theta)$
constructed in section \low~ with the condition for a locally free action
being that all the minor determinants $\grD_{ij}(\Theta)$ are nonvanishing. 

Finally we briefly discuss the two singular cases 
$$\calm(1,1,1)\simeq \bbz_3\backslash G_2/Sp(1) \quad \hbox{and} \quad
\calm(1,1,1,1)\simeq \bbz_2\backslash \hbox{Spin}(7)/\hbox{Spin}(4).$$   
Since these are  biquotients of Lie groups the topology is more accessible.
In particular, their rational cohomology is that of the corresponding
3-Sasakian homogeneous space, $G_2/Sp(1)$ and  $\hbox{Spin}(7)/\hbox{Spin}(4),$
respectively, which is well known [GS,BG2].  Thus, $\calm(1,1,1)$ has the
rational cohomology of $S^{11},$ whereas $\calm(1,1,1,1)$ has the rational
cohomology of $S^4\times S^{11}.$ In both cases $b_2$ vanishes, and we do not
expect this in the non-singular cases. 

\vfill\eject
\bigskip 
\centerline {\bf \hcx. Hypercomplex Structures on Circle Bundles over 
$\calm(\bfp)$} 
\medskip 
According to the general theory described in [BGM3] 3-Sasakian manifolds
(orbifolds) give rise to hypercomplex structures on circle bundles over them.
In this short section we give new hypercomplex structures in dimensions 12 and
16 constructed as circle V-bundles over the 3-Sasakian orbifolds $\calm(\bfp)$
constructed in sections 2 and 4. Of course, there is the trivial bundle
$\calm(\bfp)\times S^1$ over $\calm(\bfp)$ which always admits a locally
conformally hyperk\"ahler structure, but here we concentrate on the level sets
$\calh(\bfp).$ As discussed in Remark \top.8 these level sets will be smooth
manifolds as long as $0<p_1<p_2<p_3$ in the 12 dimensional case and $0\leq
p_1<p_2<p_3<p_4$ in the 16 dimensional case. We now have from our previous
results [BGM3]:

\noindent{\sc Theorem \hcx.1}: \tensl Let $\bfp$ have components satisfying the
inequalities above, then ${\calh}(\bfp)$ is a compact hypercomplex manifold of
dimension 12 or 16. Furthermore, the connected component of the Lie group of
hypercomplex automorphisms is $T^3$ in the 12 dimensional case and $T^4$ in
the 16 dimensional case. \tenrm    

The last statement of Theorem \hcx.1 implies that these hypercomplex
structures are distinct from any of those known previously. We do not know
whether for different $\bfp$ the manifolds $\calh(\bfp)$ are diffeomorphic or
not; however, arguments similar to those in [BGM2] show that the hypercomplex
structures are distinct. In fact, each smooth manifold $\calh(\bfp)$ has a real
one parameter family of distinct hypercomplex structures on them given by
sending $\bfp\mapsto \grl\bfp$ for any real $\grl>0.$

We can also construct hypercomplex structures on the total
space $\calh(\Theta)$ of circle V-bundles over the 7-dimensional 3-Sasakian
orbifolds  $\calm(\Theta)$ constructed in Section 3. In this case as in [BGM3]
there should be gcd conditions on the entries of the matrix $\Theta$ that
gaurentee that $\calh(\Theta)$ be a smooth manifold. These then give new
hypercomplex manifolds in dimension 8 with a two-dimensional group of
hypercomplex automorphisms.

\bigskip
\bigskip
\medskip
\centerline{\bf Bibliography}
\medskip
\font\ninesl=cmsl9
\font\bsc=cmcsc10 at 10truept
\parskip=1.5truept
\baselineskip=11truept
\ninerm 
\item{[BCR]} {\bsc J. Bochnak, M. Coste, and M.-F. Roy}, {\ninesl G\'eom\'etrie
alg\'ebrique r\'eele}, Springer-Verlag, Berlin, 1987.
\item{[BG1]} {\bsc C.P. Boyer and  K. Galicki}, {\ninesl
On Sasakian-Einstein geometry}, mathDG/9811098 to appear in Int. J. Math. 
\item{[BG2]} {\bsc C.P. Boyer and  K. Galicki}, {\ninesl
3-Sasakian Manifolds}, Surveys in Differential Geometry, Volume VI,
{\it Essays on Einstein Manifolds}, A supplement to
the Journal of Differential Geometry, pp. 123-184,
C. LeBrun and M. Wang, Eds., International Press, Cambridge 1999.
\item{[BGM1]} {\bsc C.P. Boyer, K. Galicki, and B.M. Mann}, {\ninesl
The geometry and topology of 3-Sasakian manifolds}, 
J. reine angew. Math., 455 (1994), 183-220.
\item{[BGM2]} {\bsc C.P. Boyer, K. Galicki, and B.M. Mann}, {\ninesl
Hypercomplex structures on Stiefel manifolds}, Ann. Global Anal. Geom. 14
(1996), 81-105.
\item{[BGM3]} {\bsc C.P. Boyer, K. Galicki, and B.M. Mann}, {\ninesl 
Hypercomplex structures from 3-Sasakian structures}, 
J. reine angew. Math., 501 (1998), 115-141.
\item{[BGMR]} {\bsc C.P. Boyer, K. Galicki, B.M. Mann, and E. Rees},
{\ninesl Compact 3-Sasakian 7-Manifolds with Arbitrary Second Betti Number},
Invent. Math. 131 (1998), 321-344.
\item{[BGP]} {\bsc C.P. Boyer, K. Galicki, and P. Piccinni},
{\ninesl Torus actions on $Gr_4(\bbr^n)$ and 3-Sasakian manifolds}, 
in preparation.
\item{[BGOP]} {\bsc C.P. Boyer, K. Galicki, L. Ornea, and P. Piccinni},
{\ninesl Geometry of Exceptional Quotients}, in preparation.
\item{[DS]} {\bsc A. Dancer and A. Swann}, {\ninesl
The Geometry of Singular Quaternionic K\"ahler
Quotients}, Int. J. Math., 8 (1997), 595-610.
\item{[G]} {\bsc K. Galicki}, {\ninesl
A generalization of the momentum mapping construction for
quaternionic K\"ahler manifolds}, Commun. Math. Phys 108 (1987), 117-138.
\item{[GL]} {\bsc K. Galicki and B. H. Lawson, Jr.}, {\ninesl
Quaternionic Reduction and Quaternionic Orbifolds}, Math. Ann., 282 (1988),
1-21.
\item{[GS]} {\bsc K. Galicki and S. Salamon}, {\ninesl On
Betti numbers of 3-Sasakian manifolds},  Geom. Ded. 63 (1996), 45-68.
\item{[Hi1]} {\bsc N. J. Hitchin}, {\ninesl A
new family of Einstein metrics.} {\it  Manifolds and geometry 
(Pisa, 1993)}, 190--222, Sympos. Math., XXXVI, 
Cambridge Univ. Press, Cambridge, 1996.
\item{[Hi2]} {\bsc N. J. Hitchin}, {\ninesl
Twistor spaces, Einstein metrics and 
isomonodromic deformations}, J. Differential Geom. 42 (1995), 30--112.
\item{[HL]} {\bsc R. Harvey and B. H. Lawson, Jr.}, {\ninesl
Calibrated geometries}, Acta Math. 148 (1982), 47--157. 
\item{[Kr1]} {\bsc  P. Kronheimer}, {\ninesl
A hyperk\"ahler structure on coadjoint orbits of a semisimple
complex group}, J. Lond. Math. Soc., II. Ser. 42 (1990), 193-208.
\item{[Kr2]} {\bsc  P. Kronheimer}, {\ninesl 
Instantons and geometry of the nilpotent variety},
J. Differential Geom. 32 (1990), 473-490.
\item{[KS1]} {\bsc  P. Kobak and A. Swann}, {\ninesl 
Quaternionic geometry of a nilpotent variety}, Math. Ann. 297 (1993),
747-764.
\item{[KS2]} {\bsc  P. Kobak and A. Swann}, {\ninesl  
Classical nilpotent orbits as hyperk\"ahler quotients},
Internat. J. Math. 7 (1996), 193-210.
\item{[KS3]} {\bsc  P. Kobak and A. Swann}, {\ninesl
Exceptional hyperk\"ahler reductions}, Twistor Newsletter 44 (1998), 23-26.
\item{[KS4]} {\bsc  P. Kobak and A. Swann}, {\ninesl Hyperk\"ahler
Potentials in Cohomogeneity Two}, math.DG/0001024.
\item{[OP]} {\bsc L. Ornea and P. Piccinni}, {\ninesl On some moment maps and
induced Hopf bundles on the quaternionic projective plane}, mathDG/0001066, to
appear in Int. J. Math.
\item{[PS]} {\bsc  Y. S. Poon and S. Salamon}, {\ninesl  
Eight-dimensional quaternionic K\"ahler manifolds with positive 
scalar curvature}, J. Differential Geom. 33 (1990), 363-378. 
\item{[Sw]} {\bsc A. F. Swann}, {\ninesl Hyperk\"ahler and
quaternionic K\"ahler geometry}, Math. Ann., 289 (1991), 421-450.
\medskip
\bigskip \line{ Department of Mathematics and Statistics
\hfil July 2000} \line{ University of New Mexico \hfil} \line{ Albuquerque, NM
87131 \hfil } \line{ email: cboyer@math.unm.edu, galicki@math.unm.edu
\hfil}
\bigskip \line{ Universit\`a degli Studi di Roma, ``La Sapienza"
\hfil } \line{ Piazzale Aldo Moro 2 \hfil} 
\line{ I-00185 Roma, Italia
\hfil } \line{ email: piccinni@mat.uniroma1.it
\hfil}
\bye